\documentclass[twoside,draft]{article}

\usepackage{amsfonts}
\usepackage{stmaryrd}
\usepackage{amssymb}
\usepackage{euscript}
\usepackage{amsthm}
\usepackage{amsmath}
\usepackage{amscd}
\usepackage{latexsym}
\usepackage{mathrsfs}
\usepackage{graphicx}
\usepackage{color}
\usepackage{dsfont}
\usepackage{bm}

\numberwithin{equation}{section}

\newtheorem{theorem}{Theorem}[section]
\newtheorem{lemma}{Lemma}[section]
\newtheorem{proposition}{Proposition}[section]

\newtheorem{remark}{Remark}[section]

\newtheorem{definition}{Definition}[section]

\newtheorem{hyp}{Assumption}[section]
\newcommand{\bhyp }{\begin{hyp} \rm }
\newcommand{\ehyp }{\end{hyp}}

\def\proof{\noindent {\it Proof. $\, $}}
\def\endproof{\hfill $\Box$}

\newcommand{\cadlag}{RCLL}
\newcommand{\cadlags}{RCLL }

\newcommand{\wt}{\widetilde }
\newcommand{\wh}{\widehat }

\newcommand{\ovl}{\overline}

\newcommand{\etab}{\eta }
\newcommand{\wtetab}{\wt{\eta }}
\newcommand{\etasq}{\eta^2}
\newcommand{\II}{I}
\newcommand{\JJ}{J}
\newcommand{\Ku}{U}
\newcommand{\Kl}{L}
\newcommand{\ovlY}{\overline{Y}}
\newcommand{\ovg}{\overline{g}}

\def\phi{\varphi}
\newcommand{\UU }{D}

\newcommand{\cTp}{\cT_p}

\newcommand{\deltag}{\delta}
\newcommand{\bigG}{G}

\newcommand{\dimn}{n}

\newcommand{\Ss}{{\mathcal{S}}^{2}}
\newcommand{\Hs}{{\mathcal{H}}^{2}(Q)}
\newcommand{\Ls}{{\mathcal{L}}^{2}(M)}
\newcommand{\Asq}{{\mathcal{A}}^{2}}
\newcommand{\Hb}{{\mathcal{H}}^{2}_{\beta}(Q)}
\newcommand{\Lb}{{\mathcal{L}}^{2}_{\beta }(M)}
\newcommand{\Sb}{{\mathcal{S}}^{2}_{\beta }}
\newcommand{\SSs}{\mathbb{S}^2}
\newcommand{\HHs}{\mathbb{H}^2}
\newcommand{\SSb}{\mathbb{S}^2_{\beta}}
\newcommand{\HHb}{\mathbb{H}^2_{\beta}}
\newcommand{\HHA}{\mathbb{H}^2_{{\mathcal A}}}
\newcommand{\SSA}{\mathbb{S}^2_{{\mathcal A}}}
\newcommand{\HHAA}{\wt{\mathbb{H}}^2_{{\mathcal A}}}
\newcommand{\SSAA}{\wt{\mathbb{S}}^2_{{\mathcal A}}}

\newcommand{\Ltg}{L^2({\cG}_T)}

\newcommand\I{\mathds{1}}

\newcommand{\E}{{\mathbb E}}

\newcommand{\esssup}{\operatornamewithlimits{ess\,sup}}

\newcommand{\cT}{\mathcal T}
\newcommand{\cG}{\mathcal F}

\newcommand{\cB}{\mathcal B}
\newcommand{\cL}{\mathcal L}
\newcommand{\cS}{\mathcal S}

\newcommand{\bff}{\mathbb F}

\newcommand{\bnn}{\mathbb N}
\newcommand{\bpp}{\mathbb P}
\newcommand{\brr}{\mathbb R}

\title{{\Large \bf REFLECTED BSDEs AND DOUBLY REFLECTED BSDEs \\ DRIVEN BY RCLL MARTINGALES} \vskip 40 pt }

\author{Tianyang Nie$\,^{a}$\footnote{The research of T. Nie and M. Rutkowski was supported by the DVC Research Bridging Support Grant {\it Pricing of American and game options in market with frictions}. The work of T. Nie was supported by the National Natural Science Foundation of China (No. 12022108,11971267,11831010,61977043).} \ and Marek Rutkowski$\,^{b,c}$ \\ \\
\\$^{a\,}$School of Mathematics, Shandong University,\\ Jinan, Shandong
250100, China\\ \\  $^{b\,}$School of Mathematics and Statistics, University of Sydney
\\ Sydney, NSW 2006, Australia\\ \\ $^{c\,}$Faculty of Mathematics and Information Science,
Warsaw University of Technology, \\ 00-661 Warszawa, Poland \\ }

\date{\vskip 30 pt \today \vskip 25 pt}

\begin{document}

\maketitle

\begin{abstract}
We prove some new results on reflected BSDEs and doubly reflected BSDEs driven by a multi-dimensional RCLL martingale.
The goal is to develop a general multi-asset framework encompassing a wide spectrum of nonlinear financial models, including as
particular cases the setups studied by Peng and Xu \cite{PX2009} and Dumitrescu et al. \cite{DGQS2018} who dealt with BSDEs driven by a one-dimensional Brownian motion and a purely discontinuous martingale with a single jump. Our results are not covered by existing literature on reflected and doubly reflected BSDEs driven by a Brownian motion and a Poisson random measure.


\end{abstract}

\renewcommand{\thefootnote}{\fnsymbol{footnote}}
\footnotetext{\textit{{E-mail:}} {nietianyang@sdu.edu.cn (Tianyang\ NIE); marek.rutkowski@sydney.edu.au (Marek\ RUTKOWSKI).}}

\newpage

\section{Introduction}   \label{sec1}
Reflected {\it backward stochastic differential equations} (BSDEs)  were introduced in seminal papers by El Karoui et al.~\cite{EQK1997,EPAQ1997} where they were applied to solutions of optimal stopping problems and pricing of American options. They were subsequently studied by several authors who dealt with various frameworks, to mention a few: Aazizi and Ouknine~\cite{AO2016}, Baadi and Ouknine~\cite{BB2017,BB2018}, Bayraktar and Yao~\cite{BY2012}, Cr\'epey and Matoussi~\cite{CM2008}, Essaky~\cite{ES2008}, Hamad\`ene~\cite{HA2002}, Hamad\`ene and Ouknine~\cite{HO2016}, Grigorova et al.~\cite{GIOOQ2017,GIOQ2017}, Klimsiak~\cite{K2012,K2015}, Klimsiak et al.~\cite{KRS2019}  Lepeltier and Xu \cite{LX2005}, Peng and Xu \cite{PX2005}, and Quenez and Sulem~\cite{QS2014}.

Doubly reflected BSDEs were first introduced and studied by Cvitani\'c and Karatzas \cite{CK1996} who also showed that classical Dynkin games can be solved using the theory of doubly reflected BSDEs. Their studies were continued by, among others, Bayraktar and Yao~\cite{BY2015}, Cr\'epey and Matoussi~\cite{CM2008}, Dumitrescu et al.~\cite{DQS2016}, Essaky and Hassani~\cite{EH2013}, Grigorova et al. \cite{GIOQ2017}, Grigorova and Quenez~\cite{GQ2017}, Hamad\`ene and Lepeltier \cite{HL2000}, Hamad\`ene and Ouknine~\cite{HO2016}, Hamad\`ene and Wang~\cite{HW2009}, Kobylanski et al. \cite{KQC2014}, and Lepeltier and San Mart\'in~\cite{LSM2004}. For a survey of results on Dynkin games and their applications to the valuation and hedging of game options, the reader is referred to Kifer~\cite{KY2013}.

Our goal is to study the class of reflected BSDEs and doubly reflected BSDEs, which are not covered by the existing literatures, but appear in a natural way in applications of reflected BSDEs and doubly reflected BSDEs to optimal stopping problems and financial mathematics. This work can be seen as a direct continuation of our previous work \cite{NR2019} where a particular class of BSDEs driven by discontinuous martingales was studied. The main motivation for this research is to provide a theoretical underpinning for stochastic models of nonlinear financial markets with credit risk where discontinuous martingales appear in a natural way. In particular, our results encompass market models with a single extraneous event, usually representing some credit event,  which were previously studied by Peng and Xu \cite{PX2009,PX2010} and, more recently, by Dumitrescu et al. \cite{DGQS2018}. The interested reader is also referred to Dumitrescu et al. \cite{DQS2016,DQS2017,DQS2018} for applications of results from \cite{DGQS2018} on reflected and doubly reflected BSDEs with a single jump to the nonlinear valuation and superhedging of American and game options in financial models with the reference credit risk.

It should be stressed that in related papers \cite{DGQS2018,PX2009,PX2010}, the authors examine a special class of BSDEs driven by a one-dimensional Brownian motion and a purely discontinuous martingale with a single jump of a unit size,  and the reflected BSDEs and doubly reflected BSDEs were  studied by Dumitrescu et al. \cite{DQS2016,DQS2017,DQS2018} in the same framework. We have extended the results of \cite{DGQS2018,PX2009,PX2010} to the BSDEs driven by a fairly general class of multi-dimensional RCLL martingales in \cite{NR2019}. It should be acknowledged that several classes of BSDEs driven by a Brownian motion, and possibly also a Poisson random measure, were studied by, among others, Barles et al. \cite{BBP1995}, Cr\'epey and Matoussi \cite{CM2008}, El Karoui et al. \cite{EMN2016}, Hamad\`ene and Ouknine \cite{HO2016}, Hamad\`ene and Wang \cite{HW2009},  Jeanblanc et al. \cite{JMN2012,JMN2013}, Papapantoleon et al. \cite{PPS2018}, Quenez and Sulem \cite{QS2013}, Royer \cite{R2006} and Tang and Li \cite{TL1994}. However, the results obtained in the above-mentioned papers do not cover the class of BSDEs driven by RCLL martingales considered in \cite{NR2019}.

In this paper,  we will extend the work \cite{NR2019} to  study the reflected BSDEs and doubly reflected BSDEs driven by a fairly general martingale satisfying Assumption \ref{ass2.1} and thus results from  \cite{DQS2016,DQS2017,DQS2018} are covered by the present work but our setup is manifestly more encompassing.  Our motivation of studying this kind of general reflected BSDEs and doubly reflected BSDEs comes from the pricing of American and game options in general nonlinear financial market models, for instance, market models with the credit risk where several default times appears.

As in \cite{NR2019}, we postulate here that the predictable representation property of the driving martingale (see Assumption \ref{ass2.2}), which is known to hold in some circumstances (see, e.g., Kusuoka \cite{K1999} or Jeanblanc and Le Cam \cite{JC2009}, in particular it also holds in the framework of \cite{DQS2016,DQS2017,DQS2018}). We stress that this assumption can be readily relaxed in theorems yielding the existence and uniqueness of solutions, provided that an additional orthogonal martingale term is introduced in the definition of a solution to a reflected or doubly reflected BSDE, as was done in Carbone et al. \cite{CFS2008} and El Karoui and Huang \cite{ELH1997} (see Remark \ref{rem2.1} of current paper and Remark  2.1 in \cite{NR2019} for more details). For conciseness, we do not elaborate here on that fairly natural extension. 

Our results are capable of covering a much broader spectrum of stochastic models, for instance, the financial market models with the credit risk of multiple entities where several default times appear in a natural way. More generally, the results obtained in this work can be employed to derive the results of the valuation, hedging and exercising of American and game options  in any market model where the price processes of primary assets are governed by stochastic differential equations driven by discontinuous noise processes. Some applications of our results to the American and game options are presented in follow-up works by Kim et al. \cite{KNR2018a,KNR2018b}.

The paper is organized as follows. We start by presenting in Section \ref{sec2} some definitions and assumptions from \cite{CFS2008,ELH1997,NR2016,NR2019} pertaining to BSDEs driven by RCLL martingales. 
In Section \ref{sec7}, we study the reflected BSDEs driven by RCLL matingales and we obtain the {\it a priori} estimates (Propositions \ref{xpro3.1}--\ref{xpro3.3}), the existence and uniqueness result (Theorem \ref{the3.1}), as well as the convergence of Picard's iterations for solutions to reflected BSDEs (Proposition \ref{xpro3.4}). Analogous results for doubly reflected BSDEs driven by RCLL martingales are established in Section \ref{sec4x}.

\section{Preliminaries}     \label{sec2}

We start by describing the setup and assumptions introduced in  \cite{NR2019} where BSDEs driven by \cadlags martingales are studied. Let $(\Omega,\cG,\bff ,\bpp)$ be a filtered probability space satisfying the usual conditions of right-continuity and completeness. The initial $\sigma$-field $\cG_0$ is assumed to be trivial.  For any two processes, say $X$ and $Y$, we denote the equality $X=Y$ as $\bpp (X_t=Y_t,\,\forall\,t\in [0,T])=1$. 
For any process $X$, we set $X_{0-}=0$ by convention, and $\Delta X_0 := X_0-X_{0-}=X_0$ is the jump of $X$ at time 0. Moreover, the integral $\int_s^t X_u\,dY_u$ is interpreted as $\int_{]s,t]}X_u\,dY_u$ where $]s,t]=\{u\in [0,T]:s<u\le t\}$.

We assume that $M=((M^1_t,M^2_t,\ldots,M^{\dimn}_t)^{\ast},\,t\in [0,T])$ is an $\dimn$-dimensional, square-integrable martingale on $(\Omega,\cG,\bff,\bpp)$. Since the filtration $\bff$ is right-continuous, an \cadlags modification of any $\bff$-martingale is known to exist and thus it is assumed throughout that $M$ is an \cadlags process. We denote by $\langle M \rangle$ (resp. $[M]$) the predictable covariation process (resp. the covariation process ) of $M$.  As in \cite{CFS2008,ELH1997,NR2016,NR2019}, we work under the following assumption regarding the martingale $M$.

\bhyp \label{ass2.1}
There exists an $\brr^{\dimn \times \dimn}$-valued, $\bff$-predictable process $m$ and an $\bff$-adapted, continuous, nondecreasing process $Q$ with $Q_0=0$ such that, for all $t\in [0,T]$,
\begin{equation} \label{eq2.1}
\langle M\rangle_t=\int_0^t m_u m_u^{\ast}\,dQ_u.
\end{equation}
In addition, we assume that $Q$ is bounded i.e. $\exists$ a constant $C_Q$ s.t. $0\le Q_t\le C_Q$ for all $t\in[0,T]$.
\ehyp
\begin{remark}
{\rm Assumption \ref{ass2.1} holds naturally under some suitable assumptions, see \cite{NR2019}  for more details. Without loss of generality, we may assume that the process $m$ appearing in \eqref{eq2.1} takes values in the space of symmetric matrices so that $m_u=m^{\ast}_u=(m_u m_u^{\ast})^{\frac{1}{2}}$,  see e.g. \cite{NR2016,NR2019}.} 


\end{remark}

Now we introduce some spaces which will be used in the sequel of the paper. As usual, $\Ltg$ stands  for the space of all real-valued, $\cG_{T}$-measurable random variables $\eta$ such that $\|\eta\|_{\Ltg}^2:=\E (\eta^2)<\infty$.
Let $\Ls$ be the space of all $\brr^{\dimn}$-valued, $\bff$-predictable processes $Z$ with the pseudo-norm
\begin{align*}
\|Z\|_{\Ls}^2:=\E\bigg[\int_0^T\|m_tZ_t\|^2\,dQ_t\bigg]<\infty,
\end{align*}
where $\|\cdot\|$ denotes the Euclidean norm in $\brr^{\dimn}$. For any fixed $\beta \ge 0$, we denote by $\Lb$ the class of all
$\brr^{\dimn}$-valued, $\bff$-predictable processes $Z$ with the pseudo-norm
\begin{align*}
\|Z\|_{\Lb}^2:=\E \bigg[\int_0^T e^{\beta Q_t}\| m_tZ_t\|^2\,dQ_t\bigg]<\infty.
\end{align*}
Hence $\Ls=\cL^2_{0}(M)$ and the pseudo-norms $\|\cdot\|_{\Ls}$ and $\|\cdot\|_{\Lb}$ are equivalent for every $\beta >0$ since the nondecreasing process $Q$ is nonnegative and bounded.

Let $\Ss$ stand for the space of all real-valued, \cadlag, $\bff$-adapted processes $X$ with the norm
\begin{align*}
\|X\|_{\Ss}^2:=\E\bigg[\sup_{t\in [0,T]}X^2_t\bigg]<\infty.
\end{align*}
Similarly, we denote by $\Sb$ the space  of all real-valued, \cadlag, $\bff$-adapted processes $X$ with the norm
\begin{align*}
\|X\|_{\Sb}^2:=\E \bigg[\sup_{t\in [0,T]}\big(e^{\beta Q_t} X^2_t\big)\bigg]<\infty.
\end{align*}
It is clear that $\Ss=\cS_{0}^2$ and the norms $\|\cdot\|_{\Ss}$ and $\|\cdot\|_{\Sb}$ are equivalent for every $\beta >0$.

We denote by $\Hs$ the space of all real-valued, $\bff$-progressively measurable processes $X$ with the pseudo-norm
\begin{align*}
\|X\|_{\Hs}^2:=\E\bigg[\int_0^T X^2_t\,dQ_t\bigg]<\infty.
\end{align*}
Similarly, we denote by $\Hb$ the space of all real-valued, $\bff$-progressively measurable processes $X$ with the pseudo-norm
\begin{align*}
\|X\|_{\Hb}^2:=\E \bigg[\int_0^T e^{\beta Q_t}X^2_t\,dQ_t\bigg]<\infty.
\end{align*}
We observe that $\Ss\subset\Hs$ and, for brevity, we denote by $\SSs$ the product space $\Ss\times\Ls$ with the following pseudo-norm,
for every $(X,Z)\in\Ss\times\Ls$,
\begin{align*}
\|(X,Z)\|_{\SSs}=\|X\|_{\Ss}+\|Z\|_{\Ls}.
\end{align*}
Similarly, we introduce the space $\HHs:=\Hs\times\Ls$ endowed with the natural pseudo-norm
\begin{align*}
\|(X,Z)\|_{\HHs}=\|X\|_{\Hs}+\|Z\|_{\Ls}.
\end{align*}
For a fixed $\beta \ge 0$, we denote by $\SSb$ the Banach space $(\Sb\times\Lb,\|\cdot\|_\beta)$ with the pseudo-norm
\begin{align*}
\|(X,Z)\|^2_\beta:=\|Z\|^2_{\Sb}+\|Z\|^2_{\Lb}
\end{align*}
and we note that $\SSs=\mathbb{S}^2_{0}$. In view of Assumption \ref{ass2.1}, the set $\SSb$ can be formally identified with $\SSs$.

We introduce the following martingale representation assumption, which provides a natural underpinning for the existence of a solution of the following BSDE \eqref{BSDE1}.

\bhyp  \label{ass2.2}
The martingale $M$ has the predictable representation property under $\bpp$ with respect to the filtration $\bff$, i.e. for any real-valued, square-integrable $\bff$-martingale $N$ with $N_0=0$ there exists a process $Z \in \Ls$ s.t. $N_t=\int_0^t Z^*_u\,dM_u$. Moreover, the uniqueness of $Z$ in $\|\cdot\|_{\Ls}$ holds, that is, if $N_t=\int_0^t Z^*_u\,dM_u=\int_0^t \wt{Z}^*_u\,dM_u$, then $\|Z-\wt{Z}\|_{\Ls}=0$.
\ehyp

\begin{remark} \label{rem2.1}
{\rm Assumption \ref{ass2.2}  is easily to be satisfied, the readers are referred to \cite{NR2019} for more comments on Assumption \ref{ass2.2}. We just mention here that, similar to Remark 2.1 in \cite{NR2019}, Assumption \ref{ass2.2} can be relaxed to the case that the generalized martingale representation property of $M$ is valid, in the sense that any real-valued, square-integrable $\bff$-martingale $N$ can be represented as $N_t=\int_0^t Z^*_u\,dM_u+L_t$ where $L$ is a square-integrable, $\bff$-martingale orthogonal to $M$.  Then the results of the existence and uniqueness of solutions in this paper still hold, provided that an additional orthogonal martingale term is introduced in the definition of a solution to a reflected or doubly reflected BSDE.} 
\end{remark}

For a fixed horizon date $T>0$, we consider the following BSDE on $[0,T]$ with the data $(g,\eta,\UU)$
\begin{equation} \label{BSDE1}
\left\{ \begin{array} [c]{ll}
dY_t=-g(t,Y_t,Z_t)\,dQ_t+Z_t^{\ast}\,dM_t+d\UU_t,\medskip\\
Y_{T}=\eta ,
\end{array} \right.
\end{equation}
where $\eta $ is a given the random variable in $\Ltg$, and  $\UU$ is a given processes in $\Hs$ s.t. $\UU_{T}\in \Ltg$. Equation \eqref{BSDE1} can be written more explicitly in the following manner, for every $t \in [0,T]$,
\begin{align} \label{BSDE1x}
Y_t=\eta+\int_t^T g(t,Y_u,Z_u)\,dQ_u-\int_t^T Z_u^{\ast}\,dM_u-(\UU_T-\UU_t).
\end{align}
Note that here $D$ is not necessary to be of finite variation.
The next definition states the measurability properties of the {\it generator} $g$ in \eqref{BSDE1}.

\begin{definition} \label{def2.1}
{\rm A {\it generator} $g:\Omega \times[0,T] \times \brr \times \brr^{\dimn } \rightarrow \brr$ is a $\cG \otimes \cB ([0,T])\otimes\cB(\brr)\otimes\cB(\brr^{\dimn})$-measurable function such that the process $g(\cdot,\cdot,y,z)$ is $\bff$-progressively measurable, for every $(y,z)\in\brr \times \brr^{\dimn}$.}
\end{definition}

Let us define $\mu_Q(A):=\E\big[\int_0^T\I_A (\omega,t)\,dQ_t(\omega)\big]$ for every $A\in\cG_T\otimes\cB([0,T])$, which is the Dol\'eans measure of the process $Q$ on the space $\Omega \times [0,T]$ endowed with the product $\sigma$-algebra $\cG_T\otimes\cB([0,T])$. We first recall the uniform Lipschitz condition for the generator, which is frequently employed in  the  study of classical BSDEs, reflected BSDEs and doubly reflected BSDEs driven by a Brownian motion (also with jumps generated by a Poisson random measure).

\begin{definition} \label{def2.2}
{\rm A generator $g$ is {\it uniformly Lipschitz continuous} if there exists a constant $L$ such that,
for $\mu_Q$-almost every $(\omega,t)$, for all $y_1,y_2\in \brr,\, z_1,z_2\in\brr^{\dimn}$,}
\begin{equation}  \label{eq2.4}
|g(t,y_1,z_1)-g(t,y_2,z_2)|\le L\big(|y_1-y_2|+\|z_1-z_2\|\big).
\end{equation}
\end{definition}

We emphasis that condition \eqref{eq2.4} is no longer suitable when dealing with BSDEs, reflected BSDEs and doubly reflected BSDEs  with a single jump, which may be used to represent the occurrence of some extraneous event $E$, and thus condition \eqref{eq2.4} has been modified to explicitly account for the {\it intensity process} $\lambda $ of $E$. For instance, in \cite{DGQS2018,DQS2017,DQS2018,PX2009,PX2010}, the authors postulate that a generator satisfies the following condition
\begin{equation}  \label{eq2.5}
|g(t,y_1,z_1,k_1)-g(t,y_2,z_2,k_2)|\le L\big(|y_1-y_2|+\|z_1-z_2\|+\sqrt{\lambda_t}|k_1-k_2|\big)
\end{equation}
where the variables $k_1,k_2$ correspond to the compensated martingale of the indicator process of the event $E$. 
More generally, in existing papers on BSDEs driven by a multi-dimensional martingale (see, e.g., \cite{CFS2008,ELH1997,NR2016}), a generator is typically postulated to satisfy some kind of the $m$-Lipschitz condition with $m$ being the process from Assumption \ref{ass2.1}. 
We refer the reader to \cite{NR2016} for further references and comments on alternative Lipschitz-type conditions, which can be used when dealing with some particular classes of generators appearing in financial applications.
In the present work, we focus on the generator satisfying the following regularity condition.

\begin{definition} \label{def2.3}
{\rm A generator $g$ is {\it uniformly $m$-Lipschitz continuous} if there exists a constant $ \wh{L} >0$ such that, for $\mu_Q$-almost every $(\omega,t)$,
and all $y_1,y_2\in\brr,\, z_1,z_2\in\brr^{\dimn}$,}
\begin{equation} \label{eq2.6}
|g(t,y_1,z_1)-g(t,y_2,z_2)|\le \wh{L}\big(|y_1-y_2|+\|m_t(z_1-z_2)\|\big).
\end{equation}
\end{definition}

\begin{remark}  \label{rem2.2} {\rm
We note that condition \eqref{eq2.5} employed in Peng and Xu \cite{PX2009,PX2010} and Dumitrescu et al. \cite{DGQS2018,DQS2018}, can be seen as a special case of \eqref{eq2.6} (for more details, see Remark 2.2 in \cite{NR2019}). The readers are also referred to \cite{NR2019} for the existence and uniqueness theorem for solutions to the BSDE \eqref{BSDE1} with generator satisfying the uniform $m$-Lipschitz condition.}
\end{remark}

\section{Reflected BSDEs Driven by RCLL Martingales}     \label{sec7}

Our aim in this section is to establish the existence and uniqueness result for the  {\it reflected backward stochastic differential equation} (RBSDE) driven by RCLL martingales. Due to various kinds of the so-called {\it minimality conditions}, which are also referred to as the {\it Skorokhod conditions},
governing the behavior of a solution when it hits the obstacle, solutions to RBSDEs studied in various papers are not always directly related to the nonlinear optimal stopping problem associated with the obstacle process playing the role of reward process. Specifically, the component $Y$ of a solution to the RBSDE is not necessarily equal to the nonlinear Snell envelope of the obstacle process, although, obviously, this is indeed the case when the proof of the existence of a solution to the RBSDE hinges on the properties of the classical Snell envelope for the standard optimal stopping problem under the linear expectation, which was studied, in particular, by El Karoui \cite{EK1981}, Kobylanski and Quenez \cite{KQ2012}, and Maingueneau \cite{M1978}.

In view of our further applications to the study of American options in nonlinear markets (see \cite{KNR2018a}), the equivalence between the component $Y$ in a solution $(Y,Z,K)$ to the RBSDE and the value process of the corresponding nonlinear optimal stopping problem is a highly desirable feature and thus it will be ensured by a particular choice of minimality conditions in \eqref{RBSDE1}.

To be more specific, in this section we consider the following RBSDE on $[0,T]$ with data $(g, \eta, \UU, \xi )$ (which satisfies the following Assumption \ref{xass3.1})
\begin{equation} \label{RBSDE1}
\left\{ \begin{array} [c]{ll}
dY_t=-g(t,Y_t,Z_t)\,dQ_t+Z^{\ast}_t\,dM_t+d\UU_t-dK_t,\ Y_T=\etab , \medskip\\ Y_t \geq \xi_t , \ \forall\, t \in [0,T], \medskip\\
\int_0^T (Y_t-\xi_t)\,dK^c_t=0 \ \mbox{\rm and}\ \Delta K^d_\tau=\Delta K^d_\tau \I_{\{ (Y-D)_{\tau-}= (\xi-D)_{\tau-}\}}, \ \forall\, \tau \in \cTp,
\end{array} \right.
\end{equation}
where $\cTp$ is the class of $\bff$-predictable stopping times taking values in $[0,T]$ and $K$ is a nondecreasing, \cadlag,  $\bff$-predictable process with $K_0=0$. The continuous and discontinuous components of sample paths of $K$ are denoted by $K^c$ and $K^d$, respectively.
\bhyp \label{xass3.1}
The quadruplet $(g,\eta ,\UU,\xi )$ is such that: \hfill \break
(i) the generator $g$ is uniformly $m$-Lipschitz continuous and the process $g(\cdot,0,0)$ belongs to $\Hs $, \hfill \break
(ii) the process $\UU$ belongs to $\Hs$ and the random variable $\eta-\UU_T$ belongs to $\Ltg$, \hfill \break
(iii) the process $\xi-\UU$ belongs to $\Ss$ and $\xi_T \le \eta $.
\ehyp

Let $\Asq$ be the class of nondecreasing, \cadlag,  $\bff$-predictable processes such that $A_0=0$ and $\E(A^2_T) < +\infty$
and let $\HHA := \HHs \times \Asq = \Hs \times \Ls \times \Asq$ and $\SSA := \SSs \times \Asq = \Ss \times \Ls \times \Asq$.

\begin{definition} \label{xdef3.1}
{\rm A {\it solution} to the RBSDE  \eqref{RBSDE1}  with data $(g, \eta, \UU, \xi )$ is a triplet of stochastic processes  $(Y,Z,K) \in \HHA $ such that conditions
\eqref{RBSDE1} are satisfied in the following sense
\begin{align*}
\bpp \bigg( Y_t= \eta+\int_t^T g(t,Y_u,Z_u)\,dQ_u-\int_t^T Z_u^{\ast}\,dM_u-(\UU_T -\UU_t)+(K_T-K_t),\ \forall\,t \in [0,T] \bigg)=1,
\end{align*}
$\bpp ( Y_t \geq \xi_t,\,\forall\,t \in [0,T] )=1$ and the Skorokhod (minimality) conditions hold: $\int_0^T (Y_t-\xi_t)\,dK^c_t=0$ and $\Delta K^d_\tau = \Delta K^d_\tau \I_{\{ (Y-D)_{\tau-}= (\xi-D)_{\tau-}\}}$ for every $\tau \in \cTp$. We say that the {\it uniqueness of a solution to}  \eqref{RBSDE1} holds if for any two solutions $(Y,Z,K)$ and $(Y',Z',K')$ to \eqref{RBSDE1} we have $\| (Y,Z,K) -(Y',Z',K') \|_{\HHA }=0$.}
\end{definition}


We define $\wt{Y}:= Y-\UU,\,\wt{\xi} := \xi-\UU,\,\wtetab := \etab-\UU_T$ and $\wt{g}(t,y,z) := g(t,y+\UU_t,z)$
and we note that $\wt{g}$ is uniformly $m$-Lipschitz continuous and, in view of Assumptions \ref{xass3.1}(i)-(ii), the process $\wt{g}(\cdot,0,0)$ belongs to $\Hs $. We also set $\wt{Z}=Z$ and $\wt{K}=K$. Then RBSDE \eqref{RBSDE1} can be written equivalently as the transformed RBSDE with data $(\wt{g}, \wt{\eta}, 0,\wt{\xi})$
\begin{equation} \label{RBSDE2}
\left\{ \begin{array} [c]{ll}
d\wt{Y}_t=-\wt{g}(t,\wt{Y}_t,\wt{Z}_t)\,dQ_t+\wt{Z}^*_t\,dM_t-d\wt{K}_t,\ \wt{Y}_T=\wtetab ,  \medskip\\ \wt{Y}_t \geq \wt{\xi}_t, \ \forall\, t \in [0,T], \medskip\\
\int_0^T(\wt{Y}_t-\wt{\xi}_t)\,d\wt{K}^c_t=0 \ \mbox{\rm and}\  \Delta \wt{K}^d_\tau=\Delta \wt{K}^d_\tau \I_{\{ \wt{Y}_{\tau-}= \wt{\xi}_{\tau-}\}}, \ \forall\, \tau \in \cTp .
\end{array} \right.
\end{equation}
It is rather straightforward to check that the following lemma is valid.

\begin{lemma} \label{xlem3.1}
Under Assumption \ref{xass3.1}, RBSDE \eqref{RBSDE1} has a unique solution $(Y,Z,K) \in \SSA$ if and only if RBSDE \eqref{RBSDE2}
has a unique solution $(\wt{Y},\wt{Z},\wt{K}) \in \SSA$ and the relationships $Y=\wt{Y}+\UU,\,Z=\wt{Z},\,K = \wt{K}$ hold.
\end{lemma}

Without loss of generality, we henceforth assume that $\UU=0$, which means that we will consider the following RBSDE
on $[0,T]$ with data $(g, \eta,0, \xi )$
\begin{equation} \label{RBSDE3}
\left\{ \begin{array} [c]{ll}
dY_t=-g(t,Y_t,Z_t)\,dQ_t+Z^{\ast}_t\,dM_t-dK_t,\ Y_T=\etab , \medskip\\  Y_t \geq \xi_t, \ \forall\, t \in [0,T],\medskip\\
\int_0^T (Y_t-\xi_t)\,dK^c_t=0 \ \mbox{\rm and}\  \Delta K^d_\tau = \Delta K^d_\tau \I_{\{Y_{\tau-}= \xi_{\tau-}\}}, \ \forall\, \tau \in \cTp.
\end{array} \right.
\end{equation}

\subsection{A Priori Estimates for RBSDEs}  \label{sec8}

We will first obtain some {\it a priori} estimates for solutions to the RBSDE.
Recall that we work under the natural assumption that the obstacle process $\xi$ satisfies $\xi_T \leq \eta$.
If we denote $\xi^+_t = \max ( \xi_t ,0)$, then it follows easily from Assumption \ref{xass3.1}(iii) that the process
$\xi^+$ belongs to $\Ss$. The following result can be seen as a generalization of Proposition 3.5 in El Karoui et al. \cite{EQK1997}.

\begin{proposition} \label{xpro3.1}
Let Assumptions \ref{ass2.1}, \ref{ass2.2} and \ref{xass3.1} be valid and let $(Y,Z,K)$ be a solution to the RBSDE  \eqref{RBSDE3}.
Then there exists a constant $C\geq 0$ such that
\begin{align*}
\|Y\|_{\Ss}^2+\|Y\|_{\Hs}^2+\|Z\|_{\Ls}^2+\|K\|_{\Ss}^2\leq C\big(\|\eta\|_{\Ltg}^{2}+\|\xi^+\|_{\Ss}^2+\|g(\cdot,0,0)\|_{\Hs}^2 \big).
\end{align*}
\end{proposition}

\proof
The It\^o formula yields
\begin{align*}
dY^2_t=2Y_{t-}\,dY_t+d[Y]_t=2Y_{t-}\,dY_t+dN_t+d\langle Y \rangle_t
\end{align*}
where the process $N_t:= [Y]_t-\langle Y \rangle_t=\sum_{0\leq s\leq t} |\Delta Y_u|^2-\langle Y^d \rangle_t$ is an $\bff$-martingale. Thus
\begin{align} \label{xeq3.4}
de^{\beta Q_t}Y^2_t&=\beta e^{\beta Q_t}Y^2_t\,dQ_t+e^{\beta Q_t}\,dY^2_t \nonumber  \\
&=e^{\beta Q_t}\left[\beta Y^2_t\,dQ_t-2Y_tg(t,Y_t,Z_t)\,dQ_t+2Y_{t-}Z^*_t\,dM_t-2Y_{t-}\,dK_t\right]\\
&+e^{\beta Q_t}\left[dN_t+\|m_tZ_t\|^2\,dQ_t\right] . \nonumber
\end{align}

\noindent {\it Step 1.} We will first show that
\begin{equation}  \label{xeq3.5}
\begin{array} [c]{ll}
& \displaystyle \sup_{t\in[0,T]}\E (Y^2_t)+\E\bigg[\int_0^T|Y_u|^2\,dQ_u+\int_0^T\|m_uZ_u\|^2\,dQ_u\bigg] \\ &
\displaystyle\leq C\,\E\bigg[\etasq+\int_0^T h_u\,dQ_u+2\int_0^T\xi^+_{u-}\,dK_u\bigg]
\end{array}
\end{equation}
where we denote $h_t:=(g(t,0,0))^2$. By integrating \eqref{xeq3.4} from $t$ to $T$, taking the expectation and using the equalities $\int_0^T(Y_u-\xi_u)\,dK^c_u=0$ and $(Y_{t-}-\xi_{t-})\Delta K^d_t =0$,
we obtain
\begin{align*}
&\E (e^{\beta Q_t}Y^2_t)+\E\bigg[\int_t^T\beta e^{\beta Q_u}Y_u^2\,dQ_u\bigg]+\E\bigg[\int_t^Te^{\beta Q_u}\|m_uZ_u\|^2\,dQ_u\bigg]\medskip\\
&= \E\bigg[e^{\beta Q_T}\etasq+2\int_t^T e^{\beta Q_u}Y_u g(u,Y_u,Z_u)\,dQ_u+2\int_t^T e^{\beta Q_u}\xi_{u-}\,dK_u\bigg].
\end{align*}
Since the mapping $g$ is uniformly $m$-Lipschitz continuous, we have
\begin{align*}
|2Y_ug(u,Y_u,Z_u)|\leq & 2\wh{L}|Y_u|\big(g(u,0,0)+|Y_u|+\|m_uZ_u\|\big)\leq h_u+(3\wh{L}^2+2\wh{L})Y^2_u +\frac{1}{2}\|m_uZ_u\|^2
\end{align*}
and thus
\begin{align*}
&\E (e^{\beta Q_t} Y^2_t)+\E\bigg[\int_t^T\beta e^{\beta Q_u}Y_u^2\,dQ_u\bigg]+\frac{1}{2}\,\E\bigg[\int_t^Te^{\beta Q_u}\|m_uZ_u\|^2\,dQ_u\bigg]\\
&\leq \E\bigg[e^{\beta Q_T}\etasq+\int_t^T e^{\beta Q_u}h_u\,dQ_u+\int_t^T (3\wh{L}^2+2\wh{L})e^{\beta Q_u}Y^2_u\,dQ_u+2\int_t^T e^{\beta Q_u}\xi^+_{u-}\,dK_u\bigg].
\end{align*}
If we take $\beta=3\wh{L}^2+2\wh{L}+1$, then we obtain the following inequality, which holds for all $t \in [0,T]$,
\begin{align*}
&\E (e^{\beta Q_t} Y^2_t)+\E\bigg[\int_t^T e^{\beta Q_u}Y_u^2\,dQ_u\bigg]+\frac{1}{2}\,\E\bigg[\int_t^Te^{\beta Q_u}\|m_uZ_u\|^2\,dQ_u\bigg]\\
&\leq \E\bigg[e^{\beta Q_T}\etasq+\int_t^T e^{\beta Q_u}h_u\,dQ_u+2\int_t^T e^{\beta Q_u}\xi^+_{u-}\,dK_u\bigg].
\end{align*}
Since $0\leq Q_t \leq C_Q$ for all $t\in[0,T]$, we have that $1\leq e^{\beta Q_t}\leq e^{\beta Q_T}\leq e^{\beta C_Q}$, and thus it is easy to deduce from the inequality above that
\begin{align*}
\sup_{t\in[0,T]}\E (Y^2_t)+\E\bigg[\int_0^T Y_u^2\,dQ_u+\int_0^T\|m_uZ_u\|^2\,dQ_u\bigg]\leq e^{\beta C_Q}\,\E\bigg[\etasq+\int_0^T h_u\,dQ_u+2\int_0^T\xi^+_{u-}\,dK_u\bigg],
\end{align*}
which shows that \eqref{xeq3.5} is valid for some constant $C \geq 0$. Note that the value of a constant $C$ will vary from place to place in the remainder of the proof.

\noindent {\it Step 2.} We will now demonstrate that there exists a constant $C \geq0$ such that
\begin{align}  \label{xeq3.6}
\|K\|_{\Ss}^2=\E(K^2_T)\leq CJ
\end{align}
where
\begin{align*}
J:=\E\bigg[\etasq+\sup_{t\in[0,T]}(\xi^+_t)^2+\int_0^T h_t\,dQ_t\bigg]=\|\eta\|_{\Ltg}^{2}+\|\xi^+\|_{\Ss}^2+\| g(\cdot,0,0)\|_{\Hs}^2.
\end{align*}
From the equality
\begin{align*}
K_T=Y_0-\etab-\int_0^T g(t,Y_t,Z_t)\,dQ_t+\int_0^T Z^{\ast}_t\,dM_t
\end{align*}
we obtain
\begin{align*}
K^2_T\leq C\bigg[Y^2_0+\etasq+\bigg(\int_0^T g(t,Y_t,Z_t)\,dQ_t\bigg)^2+\bigg(\int_0^TZ^{\ast}_t\,dM_t\bigg)^2 \bigg].
\end{align*}
Since the process $Q$ is nondecreasing and bounded, there exists a constant $C \geq 0$ (independent of $t$ and $\omega$) such that, for all $t \in [0,T]$,
\begin{align*} 
\bigg(\int_t^T g(u,Y_u,Z_u)\,dQ_u\bigg)^2\leq C\,\int_t^T g^2(u,Y_u,Z_u)\,dQ_u.
\end{align*}
Using also the inequality $g^2(u,Y_u,Z_u)\leq C (h_u+ Y^2_u+\|m_uZ_u\|^2)$, estimate \eqref{xeq3.5} and the inequality
\begin{align}  \label{xeq3.7}
2C\,\int_0^T \xi^+_{t-}\,dK_t\leq 2C^2\sup_{t\in[0,T]}(\xi^+_t)^2+\frac{1}{2}\,K^2_T,
\end{align}
we find that
\begin{align*}
\E (K^2_T)& \leq C\,\E\bigg[\etasq +\int_0^T h_t\,dQ_t+2\int_0^T \xi^+_{t-}\,dK_t\bigg]\\
&\leq C\,\E\bigg[\etasq+\int_0^T h_t\,dQ_t\bigg]+2C^2\,\E\Big[\sup_{t\in[0,T]}(\xi^+_t)^2\Big]+\frac{1}{2}\,\E (K^2_T),
\end{align*}
which implies that \eqref{xeq3.6} holds.  By combining \eqref{xeq3.5}, \eqref{xeq3.6} and \eqref{xeq3.7},
we also get
\begin{align}  \label{xeq3.8}
\sup_{t\in[0,T]}\E (Y^2_t)+\E\bigg[\int_0^T Y_u^2\,dQ_u+\int_0^T\|m_uZ_u\|^2\,dQ_u+K_T^2\bigg] \leq C J.
\end{align}

\noindent {\it Step 3.} In view of \eqref{xeq3.8}, to complete the proof of the proposition, it remains to show that
\begin{align}  \label{xeq3.9}
\|Y\|^2_{\Ss} \leq C J
\end{align}
for some constant $C\geq 0$. Since
\begin{align*}
Y_t=\etab+\int_t^T g(u,Y_u,Z_u)\,dQ_u-\int_t^TZ^{\ast}_u\,dM_u+K_T-K_t,
\end{align*}
we have, for every $t\in [0,T]$,
\begin{align*}
Y^2_t\leq C\bigg[\etasq+\bigg(\int_t^T g(u,Y_u,Z_u)\,dQ_u \bigg)^2+\bigg(\int_t^TZ^{\ast}_u\,dM_u\bigg)^2+K^2_T \bigg].
\end{align*}
The Burkholder-Davis-Gundy inequality with $p=2$ applied to the $\bff$-martingale $\wh{M}_t := \int_0^t Z^{\ast}_u\,dM_u$ gives
\begin{align*}
\|\wh{M}\|^2_{\Ss}\leq C\,\E\big([\wh{M}]_T\big)=C\,\E\bigg[\int_0^T\|m_uZ_u\|^2\,dQ_u\bigg]
\end{align*}
and thus
\begin{align*}
\|Y\|^2_{\Ss}\leq C\,\E\bigg[\etasq+\int_0^T \big(h_t+Y^2_t+\|m_tZ_t\|^2\big)\,dQ_t+K^2_T\bigg].
\end{align*}
Therefore, using \eqref{xeq3.8}, we obtain \eqref{xeq3.9} and thus the proof is completed.
\endproof

For simplicity,  we write $y:=Y^{1}-Y^{2},\,z:=Z^{1}-Z^{2}$ and $k:=K^1-K^2$. The next result furnishes another {\it a priori} estimate for solutions to the RBSDE \eqref{RBSDE3}.

\begin{proposition} \label{xpro3.2}
Let Assumptions \ref{ass2.1}, \ref{ass2.2} and \ref{xass3.1}(ii)-(iii) be valid. For $l=1,2$, assume that $(Y^l, Z^l, K^l)$ is a solution to the RBSDE
\eqref{RBSDE3} with generator $g^l$. If the mapping $g^1$ is uniformly $m$-Lipschitz continuous with the constant
$\wh{L}_1>0$, then for any $\gamma > 0$ and $\beta>0$
\begin{align*}
\big(\beta-2\wh{L}_1-\gamma^{-1}\big)\|y\|^2_{\Hb}+(1-2\wh{L}_1^2\gamma )\|z\|^2_{\Lb}\le 2 \gamma\|\ovg\|^2_{\Hb}
\end{align*}
where $\ovg_t:=g^1(t,Y^2_t,Z^2_t)-g^2(t,Y^2_t, Z^2_t)$ for all $t \in [0,T]$.
\end{proposition}

\proof
By taking the difference of the two RBSDEs, we obtain
\begin{equation*}
\left\{ \begin{array} [c]{ll}
dy_t=-\deltag_t\,dQ_t+z^*_t\,dM_t-dk_t,\medskip\\ y_T=0,
\end{array} \right.
\end{equation*}
where we denote $\deltag_t:=g^1(t,Y^1_t,Z^1_t)-g^2(t,Y^2_t,Z^2_t)$. The It\^o integration by parts formula gives
\begin{align*}
d(e^{\beta Q_t}y^2_t)=y^2_t \,de^{\beta Q_t}+e^{\beta Q_t}\,dy^2_t=y^2_t \beta e^{\beta Q_t}\,dQ_t+e^{\beta Q_t}\,dy^2_t
\end{align*}
where
\begin{align*}
dy^2_t=2y_{t-}\,dy_t+d[y]_t=2y_{t-}\,dy_t+dN_t+d\langle y \rangle_t
\end{align*}
where the process $N:=[y]-\langle y \rangle $ is an $\bff$-martingale, $d\langle y \rangle_t=\|m_t z_t\|^2\,dQ_t$, and
\begin{align*}
y_{t-}\,dy_t=-y_t\deltag_t\,dQ_t+y_{t-}z^*_t\,dM_t-y_{t-}\,dK^1_t+y_{t-}\,dK^2_t .
\end{align*}
By integrating from $t$ to $T$ and taking the conditional expectation with respect to $\cG_t$, we obtain
\begin{equation} \label{xeq3.10}
\begin{split}
&e^{\beta Q_t}y^2_t+\beta\,\E_t\bigg[\int_t^T e^{\beta Q_u}y^2_u\,dQ_u \bigg]+\E_t\bigg[\int_t^T e^{\beta Q_u}\,\|m_uz_u\|^2\,dQ_u \bigg]\\
&=2\,\E_t\bigg[\int_t^Te^{\beta Q_u}y_u\deltag_u\,dQ_u\bigg]+2\,\E_t\bigg[\int_t^Te^{\beta Q_u}y_{u-}\,dK^1_u-2\int_t^Te^{\beta Q_u}y_{u-}\,dK^2_u\bigg].
\end{split}
\end{equation}
Furthermore, we observe that
\begin{align*}
y_{u-}\,dK^{1,c}_u=y_u\,dK^{1,c}_u=\big[(Y^1_u-\xi_u)-(Y^2_u-\xi_u)\big]\,dK^{1,c}_u=-(Y^2_u-\xi_u)\,dK^{1,c}_u \le 0
\end{align*}
and
\begin{align*}
y_{u-}\,\Delta K^{1,d}_u=\big[(Y^1_{u-}-\xi_{u-})-(Y^2_{u-}-\xi_{u-})\big]\,\Delta K^{1,d}_u \le 0
\end{align*}
so that $y_{u-}\,dK^1_u \le 0$. Analogous arguments show that $y_u\,dK^{2,c}_u \ge 0$ and $y_{u-}\,\Delta K^{2,d}_u \ge 0$ and thus the inequality
$-y_{u-}\,dK^2_u \le 0$ is valid as well.

In the next step, we use these inequalities, the fact that $|\deltag_t|\le\wh{L}_1(|y_t|+\|m_tz_t \|)+|\ovg_t|$ and the elementary inequality
$2a(\wh{L}_1 b+c)\le\gamma^{-1}a^2+2\gamma (\wh{L}_1^2 b^2+c^2)$, which holds for arbitrary real numbers $a,b,c$ and $\gamma>0$. We deduce from
\eqref{xeq3.10} that
\begin{align*}
e^{\beta Q_t} y^2_t&+\beta\,\E_t \bigg[ \int_t^T e^{\beta Q_u} y^2_u\,dQ_u\bigg]+\E_t\bigg[\int_t^T e^{\beta Q_u}\|m_uz_u\|^2\,dQ_u\bigg] \\
& \le 2 \wh{L}_1\,\E_t\bigg[\int_t^T e^{\beta Q_u} y^2_u\,dQ_u\bigg]+
2\,\E_t \bigg[\int_t^T e^{\beta Q_u} |y_u|\big(\wh{L}_1\|m_uz_u\|+|\ovg_u|\big)\,dQ_u \bigg] \\
& \le 2 \wh{L}_1\,\E_t\bigg[\int_t^T e^{\beta Q_u}|y_u|^2\,dQ_u\bigg]+ \E_t\bigg[ \int_t^T e^{\beta Q_u}
\Big( \gamma^{-1} y^2_u+2\gamma \big(\wh{L}_1^2\|m_uz_u\|^2+|\ovg_u|^2\big)\Big)\,dQ_u\bigg].
\end{align*}
Rearranging the above equation and taking $t=0$, we obtain the desired inequality.
\endproof

For any $(w, v) \in \SSb$, we suppose that the reflected BSDE \eqref{RBSDE3} associated with generator $g_t:= g(t,w_t,v_t)$ has a unique solution $(Y^{w,v},Z^{w,v},K^{w,v})\in \SSb \times \Asq $. Then we have the following result, which will be used in the next section to establish the existence and uniqueness of a solution to the reflected BSDE \eqref{RBSDE3}.

\begin{proposition} \label{xpro3.3}
Let Assumptions \ref{ass2.1}, \ref{ass2.2} and \ref{xass3.1}(ii)-(iii) be valid. For $l=1,2$, we assume that $(Y^l, Z^l, K^l)$ is a solution to the RBSDE
\eqref{RBSDE3} with generator $g^l_t=g(t,w^l_t,v^l_t)$. If the mapping $g$ is uniformly $m$-Lipschitz continuous with the constant
$\wh{L}>0$, then for any $\beta> 0$
\begin{align*}
\|y\|^2_{\Sb}+\|z\|^2_{\Lb} \leq 1158 \beta^{-1} \wh{L}^2(C_Q+1)\big( \|w\|_{\Sb}^2+\|v\|^2_{\Lb}\big)
\end{align*}
where $y=Y^1-Y^2,\,z = Z^1-Z^2,\,w=w^1-w^2$ and $v=v^1-v^2$.
\end{proposition}

\proof
By taking the difference of the two RBSDEs, we obtain
\begin{equation*}
\left\{ \begin{array} [c]{ll}
dy_t=-\deltag_t\,dQ_t+z^*_t\,dM_t-dk_t,\medskip\\ y_T=0,
\end{array} \right.
\end{equation*}
where we denote $\deltag_t:=g(t,w^1_t,v^1_t)-g(t,w^2_t, v^2_t)$. The It\^o integration by parts formula gives
\begin{align*}
d(e^{\beta Q_t}y^2_t)= y^2_t\,de^{\beta Q_t}+e^{\beta Q_t}\,dy^2_t=y^2_t\beta e^{\beta Q_t}\,dQ_t+e^{\beta Q_t}\,dy^2_t.
\end{align*}
Since
\begin{align*}
dy^2_t=2 y_{t-}\,dy_t+d[y]_t=2y_{t-}\,dy_t+dN_t+d\langle y \rangle_t
\end{align*}
where the process
$N_t:= [y]_t-\langle y \rangle_t$
is an $\bff$-martingale and $d\langle y \rangle_t = \|m_t z_t\|^2\,dQ_t$, we have
\begin{equation} \label{xeq3.11}
\begin{split}
d(e^{\beta Q_t} y^2_t)&=e^{\beta Q_t}\big(\beta y^2_t\,dQ_t-2y_{t-}\deltag_t\,dQ_t+2y_{t-}z^*_t\,dM_t-2y_{t-}\,dk_t+d[y]_t\big)\\
&=e^{\beta Q_t}\big( \beta y^2_t\,dQ_t-2y_{t-}\deltag_t\,dQ_t+2y_{t-}z^*_t\,dM_t-2y_{t-}\,dk_t+dN_t+d\langle y \rangle_t\big).
\end{split}
\end{equation}
Similarly to Proposition \ref{xpro3.2}, by integrating from $t$ to $T$,  taking the expectation, using the properties that $\int_t^Te^{\beta Q_u}y_{u-}\,dk_u\leq 0$ and $N$ is an $\bff$-martingale, we obtain
\begin{equation*}
\begin{split}
&\E \bigg[e^{\beta Q_t}y^2_t+\beta\,\int_t^T e^{\beta Q_u}y^2_u\,dQ_u +\int_t^T e^{\beta Q_u}\,\|m_uz_u\|^2\,dQ_u\bigg]\leq 2\, \E\bigg[ \int_t^Te^{\beta Q_u}y_u\deltag_u\,dQ_u \bigg] \\ &\leq 2\wh{L}\,\E\bigg[ \int_t^Te^{\beta Q_u}|y_u|(|w_u|+\|m_uv_u\|)\,dQ_u \bigg] \\
&\leq 2 \lambda^{-1} \wh{L}^2\,\E\bigg[ \int_t^Te^{\beta Q_u}y^2_u\,dQ_u\bigg] +\lambda \,\E\bigg[ \int_t^Te^{\beta Q_u} |w_u|^2\,dQ_u\bigg]+\lambda\,\E\bigg[ \int_t^Te^{\beta Q_u}\|m_uv_u\|^2\,dQ_u \bigg].
\end{split}
\end{equation*}
By taking $\lambda=2\beta^{-1}\wh{L}^2$ and recalling that the process $Q$ is bounded by a constant $C_Q$, we get
\begin{equation} \label{xeq3.12}
\begin{split}
&\E \bigg[e^{\beta Q_t}y^2_t+\int_t^T e^{\beta Q_u}\,\|m_uz_u\|^2\,dQ_u\bigg]\leq 2\,\E\,\int_t^Te^{\beta Q_u}y_u\deltag_u\,dQ_u\\
&\leq \lambda\,\E\bigg[\int_t^Te^{\beta Q_u}|w_u|^2\,dQ_u\bigg]+\lambda\,\E\bigg[\int_t^Te^{\beta Q_u}\|m_uv_u\|^2\,dQ_u\bigg] \\
&\leq \lambda(C_Q+1)\big(\|w\|_{\Sb}^2+\|v\|^2_{\Lb}\big).
\end{split}
\end{equation}
Using \eqref{xeq3.11} and the noting that $[y]$ is nondecreasing and $\int_t^Te^{\beta Q_u}y_{u-}\,dk_u\leq 0$, we get
\begin{equation*}
\begin{split}
&e^{\beta Q_t}y^2_t+\beta\,\int_t^T e^{\beta Q_u}y^2_u\,dQ_u\leq 2\,\int_t^Te^{\beta Q_u}y_u\deltag_u\,dQ_u-2\int_t^Te^{\beta Q_u}y_{u-}z^*_u\,dM_u\\
&\leq 2 \lambda^{-1}\wh{L}^2 \int_t^Te^{\beta Q_u}y_u^2\,dQ_u+\lambda\int_t^Te^{\beta Q_u}|w_u|^2\,dQ_u+\lambda\int_t^Te^{\beta Q_u}\|m_uv_u\|^2\,dQ_u-2\int_t^Te^{\beta Q_u}y_{u-}z^*_u\,dM_u ,
\end{split}
\end{equation*}
which yields
\begin{equation} \label{xeq3.13}
\begin{split}
e^{\beta Q_t}y^2_t
&\leq \lambda\int_t^Te^{\beta Q_u}|w_u|^2\,dQ_u+\lambda\int_t^Te^{\beta Q_u}\|m_uv_u\|^2\,dQ_u-2\int_t^Te^{\beta Q_u}y_{u-}z^*_u\,dM_u\\
&\leq \lambda(C_Q+1)\left(\sup_{t\in[0,T]}\Big( e^{\beta Q_t}w_t^2\Big)+\int_0^Te^{\beta Q_t}\|m_tv_t\|^2\,dQ_t\right)-2\int_t^Te^{\beta Q_u}y_{u-}z^*_u\,dM_u .
\end{split}
\end{equation}
We denote $\widetilde{y}_u=y_ue^{\frac{\beta}{2} Q_u}$ and $\widetilde{M}_t := \int_0^t e^{\frac{\beta}{2} Q_u}z^{\ast}_u\,dM_u$ and we apply the Davis inequality (see, e.g., Theorem 10.24 in He et al. \cite{HWY1992}
or Theorem 11.5.5 in Cohen and Elliott \cite{CE2015} with $p=1$)
\begin{align*}
\E\left[\sup_{t\in[0,T]}\Big|\int_0^t\widetilde{y}_{u-}\,d\widetilde{M}_u \Big|\right]\le
2 \sqrt{6}\,\E\left[ \bigg( \int_0^T|\widetilde{y}_{t-}|^2\,d[\widetilde{M}]_t\bigg)^{1/2}\right] \le
6 \,\E\left[\bigg( \int_0^T|\widetilde{y}_{t-}|^2\,d[\widetilde{M}]_t \bigg)^{1/2}\right].
\end{align*}
Consequently, using also the elementary inequalities
\begin{align*}
\E \left[\bigg(\int_0^T|\widetilde{y}_{t-}|^2\,d[\widetilde{M}]_t \bigg)^{1/2}\right]
& \leq \E\bigg[\sup_{t\in[0,T]}|\widetilde{y}_t|\,[\widetilde{M}]^{\frac{1}{2}}_T\bigg]
\leq \frac{1}{48}\,\E\bigg[\sup_{t\in[0,T]}\widetilde{y}^2_t\bigg]+ 12\,\E\big([\widetilde{M}]_T\big)
\\ &=\frac{1}{48}\,\E\bigg[\sup_{t\in[0,T]}\widetilde{y}^2_t\bigg]+12\,\E\big(\langle \widetilde{M}\rangle_T\big),
\end{align*}
we obtain
\begin{align*}
\E\left[\sup_{t\in[0,T]} \Big|\int_0^t\widetilde{y}_{u-}\,d\widetilde{M}_u \Big|\right]\leq
\frac{1}{8}\,\E\bigg[ \sup_{t\in[0,T]} \widetilde{y}^2_t\bigg]+ 72 \,\E \big(\langle \widetilde{M}\rangle_T\big),
\end{align*}
which is equivalent to
\begin{align*}
\E\left[\sup_{t\in[0,T]} \Big|\int_0^te^{\beta Q_u}y_{u-}z^*_u\,dM_u\Big|\right]\leq
\frac{1}{8}\,\E\bigg[ \sup_{t\in[0,T]}\Big(e^{\beta Q_t}y^2_t\Big) \bigg]+ 72\,\E \bigg[\int_0^T e^{\beta Q_t}\,\|m_tz_t\|^2\,dQ_t \bigg].
\end{align*}
Thus, by using \eqref{xeq3.12}, we get
\begin{equation} \label{xeq3.14}
\begin{split}
&\E\left[\sup_{t\in[0,T]} \Big|\int_t^Te^{\beta Q_u}y_{u-}z^*_u\,dM_u\Big|\right]\leq
\frac{1}{4}\,\E\bigg[ \sup_{t\in[0,T]}\Big( e^{\beta Q_t}y_t^2 \Big) \bigg]+ 144\,\E \bigg[\int_0^T e^{\beta Q_t}\,\|m_tz_t\|^2\,dQ_t \bigg] \\
&\leq \frac{1}{4}\,\E\bigg[ \sup_{t\in[0,T]}\Big( e^{\beta Q_t}y_t^2\Big) \bigg]+144\lambda(C_Q+1)\big(\|w\|_{\Sb}^2+\|v\|^2_{\Lb}\big).
\end{split}
\end{equation}
From \eqref{xeq3.13} and \eqref{xeq3.14}, we obtain
\begin{equation*}
\begin{split}
\E\bigg[ \sup_{t\in[0,T]}\Big(e^{\beta Q_t}y^2_t\Big) \bigg]&\leq \lambda(C_Q+1)\big(\|w\|_{\Sb}^2+\|v\|^2_{\Lb}\big)
+2\,\E\bigg[\sup_{t\in[0,T]} \Big|\int_t^Te^{\beta Q_u}y_{u-}z^*_u\,dM_u\bigg]\\
&\leq 289\lambda(C_Q+1)\big(\|w\|_{\Sb}^2+\|v\|^2_{\Lb}\big)+\frac{1}{2}\,\E\bigg[ \sup_{t\in[0,T]}\Big(e^{\beta Q_t}y^2_t\Big) \bigg],
\end{split}
\end{equation*}
which gives
\begin{align} \label{xeq3.15}
\E\bigg[ \sup_{t\in[0,T]}\Big(e^{\beta Q_t}y^2_t\Big) \bigg]\leq 578\lambda(C_Q+1)\big(\|w\|_{\Sb}^2+\|v\|^2_{\Lb}\big).
\end{align}
By combining \eqref{xeq3.12} and \eqref{xeq3.15}, we find that (recall that $\lambda=2\beta^{-1} \wh{L}^2 $)
\begin{align*}
\|y\|^2_{\Sb}+\|z\|^2_{\Lb} \leq 1158 \beta^{-1} \wh{L}^2 (C_Q+1)\big(\|w\|_{\Sb}^2+\|v\|^2_{\Lb}\big),
\end{align*}
which ends the proof.
\endproof

\subsection{Existence and Uniqueness Theorem for Reflected BSDEs}   \label{sec9}

Our next goal is to study the existence and uniqueness of the solution to RBSDE driven by RCLL martingales. From Lemma \ref{xlem3.1}, we know that once we obtain the existence and uniqueness of the solution to RBSDE \eqref{RBSDE3}, we also have similar results for RBSDE \eqref{RBSDE1}. We henceforth assume that the obstacle $\xi$ in \eqref{RBSDE3} (or the modified obstacle $\xi-D$ in \eqref{RBSDE1}) belongs to the space $\Ss$. 
For the reader's convenience, we give the detailed demonstration of Lemma \ref{xlem3.2}, which is adapted from the proof of Lemma 2.5 in Quenez and Sulem \cite{QS2014} who in turn use, in particular, Proposition B.11 from Kobylanski and Quenez~\cite{KQ2012} (see also \cite{KQ2016}) where the case of the classical optimal stopping problem was studied.

\begin{lemma} \label{xlem3.2}
Let Assumptions \ref{ass2.1}, \ref{ass2.2} and \ref{xass3.1} be satisfied. If the generator $g$ does not depend on $(y,z)$ so that the process $g_t:= g(t,y,z)=g(t,0,0)$ belongs to $\Hs$, then the RBSDE  \eqref{RBSDE3} has a unique solution $(Y,Z,K) \in \SSA$.
\end{lemma}

\proof
To alleviate the notation, we assume in the proof of Lemma \ref{xlem3.2} that $\xi_T = \eta $. Otherwise, it suffices to introduce an auxiliary process
$\xi_T :=\xi_t\I_{\{t<T\}}+\eta\I_{\{t=T\}}$.

We denote by $\cT_{[\sigma, T]}$ the class of all stopping times $\sigma'$ such that $\sigma \le \sigma' \le T$ and $\E_\sigma[\,\cdot\,] = \E[\,\cdot \,|\,\cG_\sigma]$ and define the stochastic system $\ovl{Y}$ by setting, for every $\sigma \in \cT$,
\begin{align*}
\ovl{Y}(\sigma):=\esssup_{\tau \in \cT_{[\sigma, T]}}\E_\sigma \Big[\xi_\tau+\int_\sigma^\tau g_t \,dQ_t\Big]
 = \esssup_{\tau \in \cT_{[\sigma, T]}}\E_\sigma \big[\xi_\tau+ \bigG_{\tau}-\bigG_{\sigma }\big]
\end{align*}
where $\bigG_t := \int_0^t g_s \,dQ_s$. From the classical aggregation theorems for stochastic systems due to
Dellacherie and Lenglart \cite{DL1982} (see also El Karoui \cite{EK1981}, Kobylanski and Quenez \cite{KQ2012} and the references therein), there exists an \cadlag,  $\bff$-adapted process $\ovl{Y}$ such that the equality $\ovl{Y}_\sigma =\ovl{Y}(\sigma)$ holds for every $\sigma \in \cT$. Furthermore, the process $V:= \ovl{Y}+ \bigG $ is a supermartingale of class (D) so it has the Doob-Meyer decomposition
$V=V_0+\ovl{M}-D$ where $\ovl{M}$ is a square-integrable martingale and $D$ is a nondecreasing, \cadlag,  $\bff$-predictable process from $\Asq$. Let $D = D^c+D^d$ be the pathwise decomposition of $D$ into its continuous and discontinuous components.

By Assumption \ref{ass2.2}, the $\bff$-martingale $M$ has the predictable representation property and thus there exists an $\bff$-predictable process $\ovl{Z} \in \Ls $ such that $\ovl{M}_t = \int_0^t \ovl{Z}_s\,dM_s$. Moreover, the process $\ovl{Z}$ is unique in the pseudo-norm  $\|\cdot \|_{\Ls }$. By an application of Proposition B.11 in \cite{KQ2012}, we obtain the equality $\int_0^T (\ovl{Y}_t-\xi_t)\,dD^c_t = 0$.  In addition, we have, for every predictable stopping time $\tau \in \cTp $,
\begin{align*}
\Delta D^d_\tau =\Delta D^d_\tau \I_{\{\ovl{Y}_{\tau-}=\wt{\xi}_{\tau-}\}}.
\end{align*}
If we set $\ovl{K}:=D^c+D^d$, then it is easy to check that $(\ovl{Y},\ovl{Z},\ovl{K})$ is a solution to the RBSDE \eqref{RBSDE2}.
We also note that the process $\ovl{Y}$ is \cadlags  and, in view of Proposition \ref{xpro3.1}, the processes $\ovl{Y}$ and $\ovl{K}$ belong to $\Ss$. This completes the proof of the existence of a solution $(Y,Z,K)$ to the RBSDE \eqref{RBSDE3} with the desired properties.

To prove the uniqueness of a solution, we assume that $(Y',Z',K')$ is an arbitrary solution to the RBSDE \eqref{RBSDE2}. We will show that $Y'=\ovl{Y}$, up to indistinguishability of stochastic processes.  Inspecting the RBSDE \eqref{RBSDE3} with the fixed driver $g_t$, we see that
$Y'+\bigG$  is a supermartingale and thus, since $Y'\ge \xi$, we obtain, for every $\sigma\in\cT$ and $\tau\in\cT_{[\sigma,T]}$,
\begin{align*}
Y'_\sigma \ge \E_\sigma \big[Y'_\tau+\bigG_{\tau}-\bigG_{\sigma}\big] \ge \E_\sigma \big[\xi_\tau+\bigG_{\tau}-\bigG_{\sigma }\big].
\end{align*}
Since $\tau \in \cT_{[\sigma, T]}$ is arbitrary, this in turn yields
\begin{align*}
Y'_\sigma \ge \esssup_{\tau \in \cT_{[\sigma, T]}}\E_\sigma \big[\xi_\tau+\bigG_{\tau}-\bigG_{\sigma }\big]=\ovl{Y}_\sigma.
\end{align*}
We will now show that $Y'_\sigma \le \ovl{Y}_\sigma$. To this end, we fix $\varepsilon > 0$ and we set $\tau_\sigma^\varepsilon := \inf \{ t \ge \sigma \,|\,Y'_t \le \xi_t+ \varepsilon \}$ so that $Y'_t > \xi_t+\varepsilon$ on $[\sigma, \tau_\sigma^\varepsilon)$ and thus the process $K^c$ is constant on $[\sigma, \tau_\sigma^\varepsilon]$. The process $K^d$ is also manifestly constant on $[\sigma, \tau_\sigma^\varepsilon)$ and, in addition, $Y'_t > \xi_t+\varepsilon$ on $[\sigma, \tau_\sigma^\varepsilon)$ gives
\begin{align*}
Y'_{\tau_\sigma^\varepsilon-} \ge \xi_{\tau_\sigma^\varepsilon-}+\varepsilon > \xi_{\tau_\sigma^\varepsilon-},
\end{align*}
which implies that $\Delta K^d_{\tau_\sigma^\varepsilon} =0$ so that $K^d$ is constant on $[\sigma, \tau_\sigma^\varepsilon]$.
Therefore, the process $Y'+ \bigG$ is a martingale on the stochastic interval $[\sigma, \tau_\sigma^\varepsilon]$.
Furthermore, from the definition of $\tau_\sigma^\varepsilon$ and the right-continuity of $Y'$ and $\wt{\xi}$, we have $Y'_{\tau_\sigma^\varepsilon} \le \xi_{\tau_\sigma^\varepsilon}+\varepsilon$ and thus
\begin{align*}
Y'_\sigma = \E_\sigma \big[ Y'_{\tau_\sigma^\varepsilon}+ \bigG_{\tau_\sigma^\varepsilon}-\bigG_{\sigma } \big]\leq \E_\sigma \big[ \xi_{\tau_\sigma^\varepsilon}+ \bigG_{\tau_\sigma^\varepsilon}-\bigG_{\sigma } \big]+\varepsilon \le \ovl{Y}_\sigma+\varepsilon .
\end{align*}
Since $\varepsilon > 0$ was arbitrary, we conclude that $Y'_\sigma \le \ovl{Y}_\sigma$. The uniqueness of $Z$ and $K$ follows from the uniqueness of the Doob-Meyer decomposition and the postulated uniqueness in the predictable representation property of $M$ (see Assumption \ref{ass2.2}).
\endproof

We now proceed to show the existence and uniqueness result for the solution of the RBSDE \eqref{RBSDE3}.
Recall that for a fixed $\beta \geq 0$, we denote by $\SSb$ the Banach space $(\Sb \times \Ls ,\|\,\cdot \,\|_\beta)$ where the norm $\|\,\cdot\,\|_\beta$
is given by
\begin{align*}
\|(w,v)\|^2_\beta:=\|w\|^2_{\Sb}+\|v\|^2_{\Lb}.
\end{align*}
To set up the proof of Theorem \ref{the3.1}, we need to define some mappings. We define the mapping $\Psi : \SSb \rightarrow \SSA$ as follows: for $(w, v) \in \SSb$ we set $\Psi(w,v)=(Y^{w,v},Z^{w,v},K^{w,v})$ where $(Y^{w,v},Z^{w,v},K^{w,v})\in \SSA$ is the unique solution of the RBSDE \eqref{RBSDE3} associated with generator $g_t(\omega) := g(\omega,t,w_t(\omega),v_t(\omega))$.  In view of Lemma \ref{xlem3.2}, the mapping $\Psi$ is well-defined. In addition, we define the projection $\pi: \SSA \rightarrow \SSb$, which maps $(w,v,l) \in \SSA $ to $\pi (w,v,l):=(w,v) \in \SSb$. Finally, we define the mapping $S: \SSb\rightarrow \SSb$ by $S:= \pi \circ \Psi$. The following lemma shows that $S$ is a contraction on $\SSb$, provided that $\beta $ is sufficiently large.

\begin{lemma} \label{xlem3.3}
If Assumptions \ref{ass2.1}, \ref{ass2.2} and \ref{xass3.1} hold, then $S: \SSb\rightarrow \SSb$ is a contraction mapping for some $\beta \geq 0$ and thus $S$  has a unique fixed point, which is denoted as $(\wh{Y}, \wh{Z})$.
\end{lemma}

\proof
We consider two pairs $(w^1,v^1)$ and $(w^2,v^2)$ in $\SSb$ and, for brevity, we denote their respective images through $S$ by $(Y^1, Z^1)$ and $(Y^2, Z^2)$.  Using the {\it a priori} estimate from Proposition \ref{xpro3.3} and the fact that the generator $g$ is uniformly $m$-Lipschitz continuous with a constant $\wh{L}$, we obtain for arbitrary $\beta \geq 0$ and $\gamma >0$
\begin{align*}
\|Y\|^2_{\Sb}+\|Z\|^2_{\Lb}\le 1158 \beta^{-1}\wh{L}^2(C_Q+1)\big(\|w^1-w^2\|^2_{\Sb}+\|v^1-v^2\|^2_{\Lb}\big).
\end{align*}
If we set $\beta=2\times1158\wh{L}^2(C_Q+1)>0$ for $\beta \geq 0$, then we obtain
\begin{align*}
\|S(w^1,v^1)-S(w^2,v^2)\|_\beta^2=\|(Y^1,Z^1)-(Y^2, Z^2)\|_\beta^2\leq\frac{1}{2}\|(w^1, v^1)-(w^2, v^2)\|_\beta^2 .
\end{align*}
By the Banach theorem, the mapping $S: \SSb\rightarrow \SSb$ has a unique fixed point $(\wh{Y}, \wh{Z})$.
\endproof

The next result establishes the existence and uniqueness of a solution to the RBSDE \eqref{RBSDE3}. Let the mapping $T: \SSA \rightarrow \SSA$
be defined by $T = \Psi \circ \pi$ so that  $T(w,v,l)=(Y^{w,v}, Z^{w,v}, K^{w,v})$ for all $(w,v,l) \in \SSA$.

\begin{theorem} \label{the3.1}
If Assumptions \ref{ass2.1}, \ref{ass2.2} and \ref{xass3.1} are satisfied, then the reflected BSDE  \eqref{RBSDE3} has a unique solution $(Y,Z,K)$ in $\HHA$.
Moreover, the triplet $(Y,Z,K)$ belongs to $\SSA$ and the process $K$ is in $\Ss$.
\end{theorem}

\proof
From the definitions of $S$ and $T$, it is clear that, for every $(w,v,l) \in \SSA$,
\begin{align} \label{xeq3.16}
T(w,v,l)=(S(w,v),K^{w,v}).
\end{align}
From Lemma \ref{xlem3.3}, we know that the mapping $S$ has a unique fixed point $(\wh{Y}, \wh{Z}) \in \SSb$.
We claim that $(\wh{Y},\wh{Z},\wh{K})$ where $\wh{K} := K^{\wh{Y},\wh{Z}} \in \Asq$ is a unique fixed point of $T$.
The equality $T(\wh{Y},\wh{Z},\wh{K})=(\wh{Y},\wh{Z},\wh{K})$ follows from the equality $S(\wh{Y}, \wh{Z})=(\wh{Y}, \wh{Z})$
and \eqref{xeq3.16} with $l =\wh{K}$. Furthermore, in view of Lemma \ref{xlem3.2}, the process $\wh{K}$ belongs to $\Ss$, since
the triplet $(\wh{Y},\wh{Z},\wh{K})$ can also be interpreted as a unique solution to the RBSDE with a fixed generator
$g_t = g(t,\wh{Y},\wh{Z})$.

To check that the fixed point of $T$ is unique, it suffices to recall that the fixed point of $S$ is unique and use again \eqref{xeq3.16}.
We have thus shown that the RBSDE \eqref{RBSDE3} has a unique solution $(\wh{Y}, \wh{Z}, \wh{K})$ where $\wh{K} = K^{\wh{Y}, \wh{Z}}$.
It is also easy to verify that $(\wh{Y},\wh{Z},\wh{K})$ is also a unique fixed point of the mapping $T^k$ for every $k \geq 2$. To this end,
we notice that the equality $T^k(w,v,l)=(S^k(w,v),K^{S^{k-1}(w,v)})$ holds for every $k \geq 2$ and, by Lemma \ref{xlem3.3},
$S^k(\wh{Y},\wh{Z})=S^{k-1}(\wh{Y},\wh{Z})=(\wh{Y},\wh{Z})$. The second assertion is now an immediate consequence of Proposition \ref{xpro3.1}.
\endproof

\subsection{Picard's Iterations for Solutions to Reflected BSDEs}   \label{sec10}

In this section, we work under the assumptions of Theorem \ref{the3.1}. Our next goal is to construct Picard's sequence and to show that it converges to the unique solution $(Y,Z,K) \in \SSA $ to the RBSDE \eqref{RBSDE3}. For this purpose, we take the sequence $\{(Y^k,Z^k)\}_{k=1}^{\infty}$ where $(Y^k,Z^k)=S^k(w,v)$ for some $(w,v)\in \SSs $ and we define the sequence $\{K^k\}_{k=1}^{\infty}$ by postulating that,  for every $k \in \bnn $, the triplet $(Y^k,Z^k,K^k)\in \SSA $ is the unique solution to the RBSDE with generator $g^k_t := g(t, Y_t^{k-1},Z_t^{k-1})$ (see Lemma \ref{xlem3.2}).
Hence we have that, for all $k \in \bnn$,
\begin{equation*}
K^k_t = -\int_0^t g(u,Y_u^{k-1},Z_u^{k-1})\,dQ_u+\int_0^tZ^{k,\ast}_u\,dM_u-Y_t^k+Y_0^k ,
\end{equation*}
 $Y^k_t \geq \xi_t$ for all $t \in [0,T]$, and the continuous and discontinuous components $K^{k,c}$ and $K^{k,d}$ of $K^k$ satisfy
\begin{align*}
\int_0^T(Y^k_t-\xi_t)\,dK_t^{k,c}=0, \quad \Delta K_{\tau}^{k,d}=\Delta K_{\tau}^{k,d}\I_{\{ \wt{Y}^k_{\tau-}= \wt{\xi}_{\tau-}\}}.
\end{align*}

From Lemma  \ref{xlem3.3}, we obtain the convergence  $\lim\limits_{n\rightarrow\infty}\|(Y^k,Z^k)-(Y, Z)\|_\beta^2=0$
where $(Y,Z,K)$  is the unique solution to the RBSDE \eqref{RBSDE3}.  The following result establishes the convergence of Picard's iterations
$(Y^k,Z^k,K^k)_{k=1}^{\infty}$ to $(Y,Z,K)$ in the Banach space $\SSA$.

\begin{proposition} \label{xpro3.4}
If Assumptions \ref{ass2.1}, \ref{ass2.2} and \ref{xass3.1} are valid, then the sequence $(Y^k,Z^k,K^k)_{k=1}^{\infty}$ converges in $\SSA$ to the unique solution $(Y,Z,K)$ of the reflected BSDE \eqref{RBSDE3}, meaning that
\begin{align} \label{xeq3.17}
\lim_{k\rightarrow\infty} \|Y^k -Y\|_{\Ss }=0 ,\quad \lim_{k\rightarrow\infty} \| Z^k-Z \|_{\Ls }=0, \quad \lim_{k\rightarrow\infty} \|K^k -K\|_{\Ss }=0.
\end{align}
\end{proposition}

\proof
The properties $\lim_{k\rightarrow\infty} \|Y^k -Y\|_{\Ss }=0 $ and $\lim_{k\rightarrow\infty} \| Z^k-Z \|_{\Ls }=0$ are immediate from the Banach fixed point theorem
(see Lemma \ref{xlem3.3}). It thus suffices to show $\lim_{k\rightarrow\infty} \|K^k -K\|_{\Ss }=0$.  If we denote $\wt{Y}^k:= Y^k-Y$, $\wt{Z}^k:= Z^k-Z$ and $\wt{K}^k:= K^k-K$,
then we have
\begin{equation*}
\left\{ \begin{array} [c]{ll}
d\wt{Y}^k_t=-\wt{g}^{k-1}_t\,dQ_t+\wt{Z}^{k,*}_t\,dM_t-d\wt{K}^k_t,\medskip\\\wt{Y}^k_T=0,
\end{array} \right.
\end{equation*}
where $ \wt{g}^{k-1}_t := g(t,Y^{k-1}_t,Z^{k-1}_t)-g(t,Y_t, Z_t)$. We claim that that there exists a constant $C \geq 0$ such that
\begin{align} \label{xeq3.18}
\|\wt{Y}^k \|^2_{\Ss}+\|\wt{K}^k\|^2_{\Ss} \le C J_k
\end{align}
where
\begin{align*}
J_k & :=\|\wt{Y}^{k-1}\|^2_{\Hs }+\|\wt{Y}^k\|^2_{\Hs }+\|\wt{Z}^{k-1}\|^2_{\Ls }+\|\wt{Z}^k\|^2_{\Ls }
\end{align*}
with $\wt{Y}^{k-1}:= Y^{k-1}-Y$ and $\wt{Z}^{k-1}:= Z^{k-1}-Z$. Assuming that \eqref{xeq3.18} holds, from the convergence  $\lim\limits_{k\rightarrow\infty}\|(Y^k,Z^k)-(Y, Z)\|_\beta^2=0$, we deduce that $\lim_{k\rightarrow\infty}J_k=0$ and thus we see that \eqref{xeq3.17} is indeed a consequence of \eqref{xeq3.18}. It thus remains to show that \eqref{xeq3.18} is satisfied.

\noindent {\it Step 1.} We will first show that
\begin{align} \label{xeq3.19}
\|\wt{Y}^k\|^2_{\Ss } \le C J_k.
\end{align}
The It\^o formula gives
\begin{align*}
d|\wt{Y}^k_t|^2=2\wt{Y}^k_{t-}\,d\wt{Y}^k_t+d[\wt{Y}^k]_t
\end{align*}
and thus
\begin{align*}
d|\wt{Y}^k_t|^2=-2 \wt{Y}^k_t\wt{g}^{k-1}_t\,dQ_t+2\wt{Y}^k_{t-}\wt{Z}^{k,*}_t\,dM_t-2\wt{Y}^k_{t-}\,d\wt{K}^k_t+d[\wt{Y}^k]_t
\end{align*}
where the process
\begin{align*}
[\wt{Y}^k]_t=\sum_{0\leq u \leq t} |\Delta\wt{Y}^k_u|^2+\langle \wt{Y}^{k,c} \rangle_t
\end{align*}
is nondecreasing.  As in the proof of Proposition \ref{xpro3.2}, we infer from the Skorokhod conditions that  $\int_t^T  \wt{Y}^k_{u-}\,d\wt{K}^k_u \leq 0$ for every $t \in [0,T]$ and thus
\begin{align*}
|\wt{Y}^k_t|^2
\leq 2 \int_t^T\wt{Y}^k_u\wt{g}^{k-1}_u\,dQ_u-2 \int_t^T\wt{Y}^k_{u-}\wt{Z}^{k,*}_u\,dM_u,
\end{align*}
which gives
\begin{equation}  \label{xeq3.20}
\sup_{t\in[0,T]} |\wt{Y}^k_t|^2\leq 2\,\int_0^T|\wt{Y}^k_t\wt{g}^{k-1}_t|\,dQ_t+ 2 \sup_{t\in[0,T]} \,\Big|\int_0^t\wt{Y}^k_{u-}\,dM^k_u\Big|
\end{equation}
where the $\bff$-martingale $M^k$ is given by $M^k_t := \int_0^t \wt{Z}^{k,\ast}_u\,dM_u$ so that $\langle M^k\rangle_t = \int_0^t \|m_u \wt{Z}^k_u\|^2\,dQ_u$.
Using also the inequality $|\wt{g}^{k-1}_t| \le \wh{L} \big(| \wt{Y}^{k-1}_t|+\|m_t \wt{Z}^{k-1}_t \|\big)$, we deduce that there exists a positive constant $C$ (as usual,
 the value of $C$ may vary from place to place in the remainder of the proof) such that, for all $t \in [0,T]$,
\begin{align*}
|\wt{Y}^k_t \wt{g}^{k-1}_t| \le \wh{L}|\wt{Y}^k_t|\big(|\wt{Y}^{k-1}_t|+\|m_t\wt{Z}^{k-1}_t\|\big)
\leq C \big(|\wt{Y}^{k-1}_t|^2+ |\wt{Y}^{k}_t|^2 +\|m_t\wt{Z}^{k-1}_t \|^2 \big)
\end{align*}
and thus it is easy to see that
\begin{align} \label{xeq3.21}
\E\bigg[\int_0^T |\wt{Y}^k_t\wt{g}^{k-1}_t|\,dQ_t\bigg] \le C J_k.
\end{align}
Furthermore, an application of the Davis inequality gives 
\begin{align*}
\E\left[\sup_{t\in[0,T]}\Big|\int_0^t\wt{Y}^k_{u-}\,dM^k_u \Big|\right]\le
2 \sqrt{6} \,\E\left[ \bigg( \int_0^T|\wt{Y}^k_{u-}|^2\,d[M^k]_u \bigg)^{1/2}\right] \le
6 \,\E\left[ \bigg( \int_0^T|\wt{Y}^k_{u-}|^2\,d[M^k]_u \bigg)^{1/2}\right].
\end{align*}
Consequently, using also the elementary inequalities
\begin{align*}
\E \left[ \bigg( \int_0^T|\wt{Y}^k_{u-}|^2\,d[M^k]_u \bigg)^{1/2}\right]
& \leq \E\bigg[\sup_{t\in[0,T]}|\wt{Y}^k_t|\,[M^k]^{\frac{1}{2}}_T\bigg]
\leq \frac{1}{16}\,\E\Big[\sup_{t\in[0,T]}|\wt{Y}^k_t|^2\Big]+ 4\,\E\big([M^k]_T\big)
\\ &=\frac{1}{16}\,\E\Big[\sup_{t\in[0,T]}|\wt{Y}^k_t|^2\Big]+4\,\E\big(\langle M^k\rangle_T\big),
\end{align*}
we obtain
\begin{align}  \label{xeq3.22}
\E\left[\sup_{t\in[0,T]} \Big|\int_0^t\wt{Y}^k_{u-}\,dM^k_u \Big|\right]\leq
\frac{3}{8}\,\E\Big[ \sup_{t\in[0,T]}|\wt{Y}^k_t|^2\Big]+ 24 \,\E \big(\langle M^k\rangle_T\big).
\end{align}
It is now easy to deduce from  \eqref{xeq3.20}, \eqref{xeq3.21} and \eqref{xeq3.22} that there exists a constant $C$ such that
\begin{align*}
\|\wt{Y}^k_t\|^2_{\Ss }= \E\Big[\sup_{t\in[0,T]}|\wt{Y}^k_t|^2\Big]\le C\,\E\bigg[\int_0^T|\wt{Y}^k_t\wt{g}^{k-1}_t|\,dQ_t\bigg]+C\,\E\big(\langle M^k\rangle_T\big)\le C J_k
\end{align*}
and thus \eqref{xeq3.19} is valid.

\noindent {\it Step 2.} In the second step, we will show that
\begin{align} \label{xeq3.23}
\|\wt{K}^k\|^2_{\Ss}\le C J_k.
\end{align}
Using the equality
\begin{align*}
\wt{K}^k_t=\wt{Y}_0-\wt{Y}_t^k -\int_0^t \wt{g}^{k-1}_u\,dQ_u+\int_0^t\wt{Z}^{k,\ast}_u\,dM_u
 = \wt{Y}_0-\wt{Y}_t^k-\int_0^t \wt{g}^{k-1}_u\,dQ_u+M^k_t,
\end{align*}
we obtain
\begin{align*}
\sup_{t\in[0,T]}|\wt{K}^k_t|\le 2\,\sup_{t\in[0,T]}|\wt{Y}^k_t|+ \int_0^T |\wt{g}^{k-1}_t| \,dQ_t +\sup_{t\in[0,T]}|M^k_t|
\end{align*}
and thus
\begin{align*}
\|\wt{K}^k\|_{\Ss } = \E\Big[\sup_{t\in[0,T]}|\wt{K}^k_t|^2\Big] & \le C\,\E\bigg[\sup_{t\in[0,T]}|\wt{Y}^k_t|^2+\Big( \int_0^T |\wt{g}^{k-1}_t|\,dQ_t \Big)^2
+ \sup_{t\in[0,T]}|M^k_t|^2 \bigg] \\
& \le C\,\E\bigg[ J_k+\int_0^T |\wt{g}^{k-1}_t|^2\,dQ_t+\int_0^T \|m_t\wt{Z}^k_t\|^2\,dQ_t \bigg]
\end{align*}
where we used \eqref{xeq3.19}, the boundedness of $Q$ and we observed that the Burkholder-Davis-Gundy inequality with $p=2$
(see, e.g., Theorem 10.36 in He et al. \cite{HWY1992} or Theorem 11.5.5 in Cohen and Elliott \cite{CE2015})  yields (one may take here $C=4$)
\begin{align*}
\|M^k\|^2_{\Ss } = \E \Big[\,\sup_{t\in[0,T]}\big|M^k_t\big|^2 \Big] \leq C \,\E ([M^k]_T)=C\,\E (\langle M^k\rangle_T)= C\,\E\bigg[\int_0^T \|m_t\wt{Z}^k_t\|^2\,dQ_t\bigg] = \|\wt{Z}^k\|^2_{\Ls }.
\end{align*}
Using the inequality $|\wt{g}^{k-1}_t|^2\le C \big(|\wt{Y}^{k-1}_t|^2+\|m_t\wt{Z}^{k-1}_t\|^2\big)$, we obtain \eqref{xeq3.23}. To conclude that inequality \eqref{xeq3.18}
is valid, it suffices to combine \eqref{xeq3.19} with \eqref{xeq3.23}.
\endproof

\section{Doubly Reflected BSDEs with RCLL Martingales}  \label{sec4x}

The goal of this section is to study the {\it doubly reflected backward stochastic differential equation} (DRBSDE) on $[0,T]$ with data $(g,\eta,\UU,\xi,\zeta)$
\begin{equation} \label{DRBSDE1}
\left\{ \begin{array} [c]{ll}
dY_t=-g(t,Y_t,Z_t)\,dQ_t+Z^*_t\,dM_t+d\UU_t-d\Kl_t+d\Ku_t,\ Y_T=\eta, \medskip \\
\xi_t\le Y_t\le \zeta_t ,\ \forall\,t \in [0,T],\medskip \\
\int_0^T (Y_t-\xi_t)\,d\Kl^c_t=\int_0^T (\zeta_t-Y_t)\,d\Ku^c_t=0,\medskip \\
\Delta \Kl^d_\tau=\Delta \Kl^d_\tau \I_{\{ (Y-D)_{\tau-}=(\xi-D)_{\tau-}\}}\ \mbox{\rm and}
\ \Delta \Ku^d_\tau=\Delta \Ku^d_\tau \I_{\{ (Y-D)_{\tau-}=(\zeta-D)_{\tau-}\}},\ \forall\,\tau\in\cTp,
\end{array} \right.
\end{equation}
where the data $(g,\eta,\UU,\xi,\zeta)$ satisfies the following assumption

\bhyp \label{xass4.1}
The quintuplet $(g, \eta, \UU, \xi , \zeta )$ is such that: \hfill \break
(i) the generator $g$ is uniformly $m$-Lipschitz continuous and the process $g(\cdot,0,0)$ belongs to $\Hs $, \hfill \break
(ii) the process $\UU$ belongs to $\Hs$ and the random variable $\eta-\UU_T $ belongs to $\Ltg$, \hfill \break
(iii) the processes $\xi-\UU$ and $\zeta-\UU$ belong to $\Ss$ and satisfy $\xi_t-\UU_t\le \zeta_t-\UU_t$ for all $t \in [0,T]$; in addition,
$\xi_T\le\eta\le\zeta_T$.
\ehyp

In \eqref{DRBSDE1}, $\cTp$ denotes all $\bff$-predictable stopping times taking values in $[0,T]$, $\Kl$ and $\Ku$ are non-decreasing, \cadlag, $\bff$-predictable processes with $\Kl_0=\Ku_0=0$ and such that the measures generated by $L$ and $U$ satisfy $d\Kl_t \perp d\Ku_t$ where the notation $d\Kl_t \perp d\Ku_t$ means that the measures generated by nondecreasing processes $\Kl$ and $\Ku$ are {\it singular},
in the sense that there exists an $\bff$-predictable set $F \subset \Omega \times [0,T]$ such that
\begin{align*}
\int_0^T \I_{F(t)}\,d\Kl_t=\int_0^T \I_{F^c(t)}\,d\Ku_t=0
\end{align*}
where $F(t):=\{ \omega \in \Omega\,|\,(\omega ,t) \in F \}$ and $F^c(t)$ is its complement.
Note that here the equality $\Kl=\Kl^c+\Kl^d$ (respectively, $\Ku=\Ku^c+\Ku^d$) gives the pathwise decomposition of $\Kl$  (respectively, $\Ku$) into its continuous and jump components. The last two rows in \eqref{DRBSDE1} are collectively referred to as the {\it Skorokhod}  (or the {\it minimality}) conditions for DRBSDEs.

To examine solutions to the RDBSDE \eqref{DRBSDE1}, we introduce the spaces of stochastic processes $\HHAA=\HHA \times \Asq$ and $\SSAA=\SSA\times\Asq$ (recall that the spaces  $\HHA,\SSA,\Asq$  are introduced in Section \ref{sec7}).

\begin{definition} \label{xdef4.1}
{\rm A {\it solution} to the DRBSDE \eqref{DRBSDE1} with data $(g,\eta,\UU,\xi,\zeta )$ is a quadruplet of stochastic processes  $(Y,Z,\Kl,\Ku)\in\HHAA $ satisfying \eqref{DRBSDE1}, in the sense that:
\begin{align*}
\bpp \bigg(Y_t=\eta+\int_t^T g(t,Y_u,Z_u)\,dQ_u-\int_t^T Z_u^{\ast}\,dM_u-(\UU_T-\UU_t)+(\Kl_T-\Kl_t)-(\Ku_T-\Ku_t),\ \forall\,t\in [0,T]\bigg)=1,
\end{align*}
$\bpp (\xi_t\le Y_t\le\zeta_t,\,\forall\,t\in [0,T])=1$ and the Skorokhod (minimality) conditions stated in \eqref{DRBSDE1} hold.
We say that the {\it uniqueness of a solution to} \eqref{DRBSDE1} holds if for any two solutions $(Y,Z,\Kl,\Ku)$ and $(Y',Z',\Kl',\Ku')$ to \eqref{DRBSDE1} we have $\| (Y,Z,\Kl,\Ku)-(Y',Z',\Kl',\Ku')\|_{\HHAA}=0$.}
\end{definition}

It can be observed that, without loss of generality, it suffices to consider the modified obstacles
$$
\xi'_t:=\xi_t\I_{\{t<T\}}+\eta\I_{\{t=T\}},\quad\zeta'_t:=\zeta_t\I_{\{t<T\}}+\eta\I_{\{t=T\}},
$$
and thus the random variable $\eta$ can be omitted, provided that we postulate that $\xi_T=\zeta_T $ in \eqref{DRBSDE1}.

As in preceding sections, we make use of the following transformations: $\wt{Y}:=Y-\UU$, $\wt{\xi}=\xi'-\UU$,  $\wt{\zeta}=\zeta'-\UU$ and $\wt{g}(t, y,z):=g(t, y+\UU_t, z)$ and we set $\wt{Z}=Z,\,\wt{\Kl}=\Kl,\,\wt{\Ku}=\Ku$. Therefore, we henceforth study  a solution $(\wt{Y},\wt{Z},\wt{\Kl},\wt{\Ku}) \in \SSAA $ to the transformed DRBSDE with data $(\wt{g},\wt{\xi }_T,0, \wt{\xi },\wt{\zeta} )$
\begin{equation} \label{DRBSDE2}
\left\{ \begin{array} [c]{ll}
d\wt{Y}_t=-\wt{g}(t,\wt{Y}_t,\wt{Z}_t)\,dQ_t+\wt{Z}^*_t\,dM_t-d\wt{\Kl}_t+d\wt{\Ku}_t,\ \wt{Y}_T=\wt{\eta }_T,  \medskip\\
\wt{\xi}_t\le \wt{Y}_t\le \wt{\zeta}_t ,\ \forall\,t \in [0,T], \medskip\\
\int_0^T (\wt{Y}_t-\wt{\xi}_t)\,d\Kl^c_t=\int_0^T (\wt{\zeta}_t-\wt{Y}_t)\,d\wt{\Ku}^c_t=0, \medskip\\
\Delta \wt{\Kl}^d_\tau=\Delta \wt{\Kl}^d_\tau \I_{\{\wt{Y}_{\tau-}=\wt{\xi}_{\tau-}\}}
\ \mbox{\rm and}\ \Delta \wt{\Ku}^d_\tau=\Delta \wt{\Ku}^d_\tau\I_{\{\wt{Y}_{\tau-}=\wt{\zeta}_{\tau-}\}},\ \forall\,\tau\in\cTp ,
\end{array} \right.
\end{equation}
and we argue that DRBSDEs \eqref{DRBSDE1} and \eqref{DRBSDE2} are equivalent under Assumption \ref{xass4.1}. Notice that, in view of inequalities
$\wt{\xi}\le \wt{Y}\le \wt{\zeta}$ and Assumption \ref{xass4.1}(iii), it is readily seen that the process $\wt{Y}$ necessarily belongs to the space $\Ss$.

\subsection{A Priori Estimates for Doubly Reflected BSDEs}  \label{xsec4.1}

To address the case of a uniformly $m$-Lipschitz continuous generator $g$, we will prove some auxiliary result furnishing the {\it a priori} estimates for solutions to \eqref{DRBSDE2}. As we argued before, we may and do assume, without loss of generality, that $\UU=0$ and $\xi_T=\zeta_T=\eta $ and thus we henceforth examine the DRBSDE with data $(g,\xi_T,0,\xi,\zeta)$. 

\begin{proposition} \label{pro4.1}
For $i=1,2$, let $(Y^i,Z^i,\Kl^i,\Ku^i )$ be a solution to the doubly reflected BSDE with data $(g^i,\xi^i_T,0,\xi^i,\zeta^i)$ and a uniformly $m$-Lipschitz
continuous generator $g^i$. If we denote $y:=Y^1-Y^2$, $z:=Z^1-Z^2$ and $\ovl{g}_t(\omega) :=g^1(\omega,t,Y^2_t,Z^2_t)-g^2(\omega,t,Y^2_t,Z^2_t)$, then for
any constant $\gamma > 0$
\begin{align*}
(\beta-2\wh{L}_1-\gamma^{-1})\|y\|^2_{\Hb }+(1-2\wh{L}_1^2 \gamma)\|z\|^2_{\Lb}\le 2\gamma\|\ovl{g}\|^2_{\Hb }
\end{align*}
where $\wh{L}_1$ is the Lipschitz constant associated with generator $g^1$.
\end{proposition}

\proof
By taking the difference of the two DRBSDEs, we obtain
\begin{equation*}
\left\{ \begin{array} [c]{ll}
y_t=-\delta_t\,dQ_t+z^*_t\,dM_t-(d\Kl^1_t-d\Kl^2_t)+(d\Ku^1_t-d\Ku^2_t),\medskip\\ y_T=0,
\end{array} \right.
\end{equation*}
where $\delta_t:=g^1(t,Y^1_t,Z^1_t)-g^2(t,Y^2_t,Z^2_t)$.  As in the proof of Proposition \ref{xpro3.2}, we obtain
\begin{equation*}
\begin{split}
&e^{\beta Q_t}y_t^2+\beta\E_t\bigg[\int_t^T e^{\beta Q_u}y^2_u\,dQ_u\bigg]+\E_t \bigg[\int_t^T e^{\beta Q_u}\| m_uz_u \|^2\,dQ_u \bigg] \\
&=2\,\E_t\bigg[\int_t^T e^{\beta Q_u}y_u\delta_u\,dQ_u\bigg]+2 \E_t\bigg[\int_t^T e^{\beta Q_u}y_{u-}\,dL^1_u-\int_t^T e^{\beta Q_u}y_{u-}\,dL^2_u\bigg] \\
&-2\,\E_t\bigg[\int_t^T e^{\beta Q_u}y_{u-}\,dU^1_u-\int_t^Te^{\beta Q_u} y_{u-}\,dU^2_u \bigg].
\end{split}
\end{equation*}
On the one hand, using the Skorokhod conditions, we observe that
\begin{align*}
y_{u-}\,dL^{1,c}_u=y_u\,dL^{1,c}_u=\big[(Y^1_u-\xi_u)-(Y^2_u-\xi_u)\big]\,dL^{1,c}_u=-(Y^2_u-\xi_u)\,dL^{1,c}_u \le 0
\end{align*}
and
\begin{align*}
y_{u-}\,\Delta L^{1,d}_u=\big[(Y^1_{u-}-\xi_{u-})-(Y^2_{u-}-\xi_{u-})\big]\,\Delta L^{1,d}_u \le 0
\end{align*}
so that $y_{u-}\,dL^1_u \le 0$. Analogous arguments show that $y_u\,dL^{2,c}_u \ge 0$ and $y_{u-}\,\Delta L^{2,d}_u \ge 0$ and thus the inequality
$-y_{u-}\,dL^2_u \le 0$ is valid as well. On the other hand, we observe that
\begin{align*}
y_{u-}\,dU^{1,c}_u=y_u\,dU^{1,c}_u=\big[(Y^1_u-\zeta_u)-(Y^2_u-\zeta_u)\big]\,dU^{1,c}_u=-(Y^2_u-\zeta_u)\,dU^{1,c}_u \ge 0
\end{align*}
and
\begin{align*}
y_{u-}\,\Delta U^{1,d}_u=\big[(Y^1_{u-}-\zeta_{u-})-(Y^2_{u-}-\zeta_{u-})\big]\,\Delta U^{1,d}_u \ge 0
\end{align*}
so that $y_{u-}\,dU^1_u \ge 0$. Similarly, one can show that $y_u\,dU^{2,c}_u \le 0$ and $y_{u-}\,\Delta U^{2,d}_u \le 0$ and thus we conclude that
$-y_{u-}\,dU^2_u \ge 0$. Therefore, we obtain
\begin{equation*}
e^{\beta Q_t}y_t^2+\beta\E_t\bigg[\int_t^T e^{\beta Q_u}y^2_u\,dQ_u\bigg]+\E_t \bigg[\int_t^T e^{\beta Q_u}\| m_uz_u \|^2\,dQ_u \bigg] \leq 2\,\E_t\bigg[\int_t^T e^{\beta Q_u}y_u\delta_u\,dQ_u\bigg].
\end{equation*}
Then, arguing as in the proof of Proposition \ref{xpro3.2}, we obtain the desired inequality.
\endproof

The next result is a counterpart of Proposition \ref{xpro3.3}.

\begin{proposition} \label{xpro4.2}
Let Assumptions \ref{ass2.1}, \ref{ass2.2} and \ref{xass4.1}(ii)-(iii) be valid. For $l=1,2$, assume that $(Y^l, Z^l, K^l)$ is a solution to the doubly
reflected BSDE \eqref{DRBSDE1} with generator $g^l_t=g(t,w^l_t,v^l_t)$. If the mapping $g$ is uniformly $m$-Lipschitz continuous with the constant
$\wh{L}>0$, then for any $\beta> 0$
\begin{align} \label{priorix}
\|y\|^2_{\Sb}+\|z\|^2_{\Lb}\le 1158 \beta^{-1} \wh{L}^2(C_Q+1)\big( \|w\|_{\Sb}^2+\|v\|^2_{\Lb}\big)
\end{align}
where $y=Y^1-Y^2,\,z=Z^1-Z^2,\,w=w^1-w^2$ and $v=v^1-v^2$.
\end{proposition}

\proof
By taking the difference of the two DRBSDEs, we obtain
\begin{equation*}
\left\{ \begin{array} [c]{ll}
y_t=-\delta_t\,dQ_t+z^*_t\,dM_t-(d\Kl^1_t-d\Kl^2_t)+(d\Ku^1_t-d\Ku^2_t),\medskip\\ y_T=0,
\end{array} \right.
\end{equation*}
where $\delta_t:=g^1(t,w^1_t,v^1_t)-g^2(t,w^2_t,v^2_t)$.  As in the proof of Proposition \ref{xpro3.3}, the It\^o integration by parts formula gives
\begin{equation}\label{xeq3.11n}
\begin{split}
d(e^{\beta Q_t} y^2_t)&=e^{\beta Q_t}\big( \beta y^2_t\,dQ_t-2y_{t-}\delta_t\,dQ_t+2y_{t-}z^*_t\,dM_t+dN_t+d\langle y \rangle_t\big)\\
&\qquad-2y_{t-}e^{\beta Q_t}\,(d\Kl^1_t-d\Kl^2_t-d\Ku^1_t+d\Ku^2_t)\\
\end{split}
\end{equation}
where the process
$N_t:= [y]_t-\langle y \rangle_t$
is an $\bff$-martingale.
By integrating from $t$ to $T$, taking the expectation, using the property that $\int_t^Te^{\beta Q_u}y_{u-}\,(d\Kl^1_u-d\Kl^2_u-d\Ku^1_u+d\Ku^2_u)\leq 0$, we obtain
\begin{equation*}
\begin{split}
&\E \bigg[e^{\beta Q_t}y^2_t+\beta\,\int_t^T e^{\beta Q_u}y^2_u\,dQ_u +\int_t^T e^{\beta Q_u}\,\|m_uz_u\|^2\,dQ_u\bigg]\leq 2\E\bigg[ \int_t^Te^{\beta Q_u}y_u\delta_u\,dQ_u \bigg] \\ &\leq 2\wh{L}\,\E\bigg[ \int_t^Te^{\beta Q_u}|y_u|(|w_u|+\|m_uv_u\|)\,dQ_u \bigg] \\
&\leq 2 \lambda^{-1} \wh{L}^2\,\E\bigg[ \int_t^Te^{\beta Q_u}y^2_u\,dQ_u\bigg] +\lambda \,\E\bigg[ \int_t^Te^{\beta Q_u}w^2_u\,dQ_u\bigg]+\lambda\,\E\bigg[ \int_t^Te^{\beta Q_u}\|m_uv_u\|^2\,dQ_u \bigg].
\end{split}
\end{equation*}
By taking $\lambda=2\beta^{-1} \wh{L}^2$ and recalling that the process $Q$ is bounded by a constant $C_Q$, we get
\begin{equation} \label{xeq3.12n}
\begin{split}
&\E \bigg[e^{\beta Q_t}y^2_t+\int_t^T e^{\beta Q_u}\,\|m_uz_u\|^2\,dQ_u\bigg]\leq 2\,\E\,\int_t^Te^{\beta Q_u}y_u\delta_u\,dQ_u\\
&\leq \lambda\,\E\bigg[\int_t^Te^{\beta Q_u}|w_u|^2 \,dQ_u\bigg]+\lambda\,\E\bigg[\int_t^Te^{\beta Q_u}\|m_uv_u\|^2\,dQ_u\bigg] \\
&\leq \lambda(C_Q+1)\big(\|w\|_{\Sb}^2+\|v\|^2_{\Lb}\big).
\end{split}
\end{equation}
Using \eqref{xeq3.11n} and noting that $[y]$ is nondecreasing and $\int_t^Te^{\beta Q_u}y_{u-}\,(d\Kl^1_u-d\Kl^2_u-d\Ku^1_u+d\Ku^2_u)\leq 0$, we obtain
\begin{equation*}
\begin{split}
&e^{\beta Q_t}y^2_t+\beta\,\int_t^T e^{\beta Q_u}y^2_u\,dQ_u\leq 2\,\int_t^Te^{\beta Q_u}y_u\delta_u\,dQ_u-2\int_t^Te^{\beta Q_u}y_{u-}z^*_u\,dM_u\\
&\leq 2 \lambda^{-1}\wh{L}^2 \int_t^Te^{\beta Q_u}y_u^2\,dQ_u+\lambda\int_t^Te^{\beta Q_u}|w_u|^2\,dQ_u+\lambda\int_t^Te^{\beta Q_u}\|m_uv_u\|^2\,dQ_u-2\int_t^Te^{\beta Q_u}y_{u-}z^*_u\,dM_u,
\end{split}
\end{equation*}
Arguing as in the proof of Proposition \ref{xpro3.3} and using \eqref{xeq3.12n} and the Davis inequality, we can
show that \eqref{priorix} is valid.
\endproof

\subsection{Existence and Uniqueness Theorem for Doubly Reflected BSDEs}  \label{sec4.2}

Before stating the next result, we need to introduce some notations. For any process $g \in \Hs$, we define
\begin{equation} \label{wt.xi}
\wt{\xi}^g_t :=\wt{\xi}_t-\E\bigg[\wt{\xi}_T+\int_t^T g_u\,dQ_u\,\Big|\,\cG_t\bigg], \quad
\wt{\zeta}^g_t :=\wt{\zeta}_t-\E\bigg[\wt{\zeta}_T+\int_t^T g_u\,dQ_u\,\Big|\,\cG_t\bigg],
\end{equation}
so that $\wt{\xi}^g_T=\wt{\zeta}^g_T=0$. If $\wt{\zeta}$ and $\wt{\xi}$ are in $\Ss $, then the processes $\wt{\xi}^g$ and $\wt{\zeta}^g$ belong to $\Ss$ as well. We define the sequences $\{\II^{g,k}\}_{k=0}^{\infty} $ and $\{\JJ^{g,k}\}_{k=0}^{\infty}$ by setting $\II^{g,0}=\JJ^{g, 0}=0$ and, by recursion, for every stopping time $\sigma \in \cT $,
\begin{align*}
\II^{g,k+1}_\sigma :=\esssup_{\tau \in \cT_{[\sigma, T]} } \E \big[ \JJ^{g,k}_\tau+\wt{\xi}^g_\tau\,|\,\cG_\sigma \big],
\quad \JJ^{g,k+1}_\sigma :=\esssup_{\tau \in \cT_{[\sigma, T]} } \E \big[ \II^{g,k}_\tau-\wt{\zeta}^g_\tau\,|\,\cG_\sigma \big]
\end{align*}
where $\II^{g,k+1}$ and $\JJ^{g,k+1}$ are taken to be RCLL. The following lemma was proven in Dumitrescu et al. \cite{DQS2016} in the special case where $Q_t=t$.  The following minor extension of Lemma 3.3 in \cite{DQS2016} can be established by applying
to the processes $\wt{\xi}^g_t$ and $\wt{\zeta}^g_t$ given by \eqref{wt.xi} the method of the proof of Lemma 3.3 in \cite{DQS2016}.

\begin{lemma} \label{ylem4.1}
The sequences of processes $\{\II^{g,k}\}_{k=0}^{\infty} $ and $\{\JJ^{g,k}\}_{k=0}^{\infty}$ are non-decreasing. Moreover, the processes $\II^g$ and $\JJ^g$, which are defined by, for all $t \in [0,T]$,
\begin{align*}
\II^g_t :=\lim_{k \rightarrow \infty}\II^{g,k},\quad \JJ^g_t :=\lim_{k \rightarrow \infty} \JJ^{g,k}
\end{align*}
are $\brr \cup \{+\infty\}$-valued, strong supermartingales with $\II^g_T=\JJ^g_T=0$ and, for every $\sigma \in \cT$,
\begin{equation} \label{Jg}
\II^g_\sigma=\esssup_{\tau\in\cT_{[\sigma, T]}}\E\big[\JJ^g_\tau+\wt{\xi}^g_\tau\,|\,\cG_\sigma\big],
\quad \JJ^g_\sigma=\esssup_{\tau \in \cT_{[\sigma, T]}}\E \big[\II^g_\tau-\wt{\zeta}^g_\tau\,|\,\cG_\sigma \big].
\end{equation}
If $\II^g_0< +\infty$ and $\JJ^g_0 < +\infty$, then $\II^g$ and $\JJ^g$ are real-valued, \cadlags supermartingales.
\end{lemma}

The proof of Lemma 3.3 in \cite{DQS2016} uses also some results on a general optimal stopping problem, which was studied, in particular, by Maingueneau \cite{M1978}, El Karoui \cite{EK1981}, Kobylanski and Quenez \cite{KQ2012} and Kobylanski et al. \cite{KQC2014}. We observe that some of their arguments can be applied to our setup. Let $\phi$ be an \cadlags process of class (D) and let the \cadlags $\bff$-supermartingale $v$ be obtained through the {\it aggregation} of the supermartingale family $\bar{v}:\cT\to\brr$, which is given by
\begin{align*}
\bar{v} (\sigma):=\esssup_{\tau\in\cT_{[\sigma,T]}}\E[\phi_\tau\,|\,\cG_\sigma].
\end{align*}
Let $B$ be the $\bff$-predictable, non-decreasing process in the Doob-Meyer decomposition of $v$. If $B=B^c+B^d$ is the pathwise decomposition of $B$ into its continuous and jump components, then it can be shown that (see, for instance, Proposition B11 in \cite{KQ2012})
\begin{align}  \label{skor}
\int_0^T \I_{\{v_t>\phi_t\}}\,dB^c_t=0\ \ \mbox{\rm and}\ \ \Delta B^d_t=\Delta B^d_t\I_{\{v_{t-}=\phi_{t-}\}}.
\end{align}
Moreover, if the filtration $\bff$ is quasi-left-continuous, then $\Delta v_\tau=\Delta B^d_\tau \I_{\{v_{\tau-}=\phi_{\tau-}\}}$ for any jump time $\tau\in\cTp $ of $v$.

We first consider the case where the generator $g$ in \eqref{DRBSDE2} is independent of $(y,z)$ and the process $g: \Omega \times [0,T] \rightarrow \brr $ belongs to $\Hs $. To show the existence and uniqueness of a solution to the RBSDE \eqref{DRBSDE2} under the assumption that  $\II^g$ and $\JJ^g$ belong to $\Ss$, we may directly apply the method of the proof of Theorem 3.5 in Dumitrescu et al. \cite{DQS2016} and thus the proof of Lemma \ref{lem4.2} is omitted.

\begin{lemma} \label{lem4.2}
Let Assumptions \ref{ass2.1}--\ref{ass2.2} and \ref{xass4.1} be satisfied and the generator $g$ do not depend on $(y,z)$. Suppose that the processes $\II^g$ and $\JJ^g$ belong to $\Ss$ and let the \cadlags process $\ovlY$ be given by, for all $t \in [0,T]$,
\begin{align} \label{Ybar}
\ovlY_t :=\II^g_t-\JJ^g_t+\E\bigg[\xi_T+\int_t^T g_u\,dQ_u\,\Big|\,\cG_t\bigg].
\end{align}
Then there exists a unique solution $(\ovlY,Z,\Ku,\Kl) \in \SSAA $ to the doubly reflected  BSDE \eqref{DRBSDE2}.
\end{lemma}

Lemma \ref{lem4.2} is not completely satisfactory, since it requires that the processes $\II^g$ and $\JJ^g$ belong to $\Ss$.
However, it is known that the following condition, which was used in several papers on the theory of Dynkin games and DRBSDEs,
ensures the existence of the value of a Dynkin game and a solution to a DRBSDE.

\bhyp \label{ass4.2} {\rm
Assume that the processes $\wt{\xi}$ and $\wt{\zeta}$ belong to $\Ss$, the inequality $\wt{\xi}_t\le \wt{\zeta}_t$ holds for all $t \in [0,T]$ and $\wt{\xi}_T= \wt{\zeta}_T$. Let the \emph{Mokobodzki condition} holds, that is, there exist two nonnegative supermartingales $H$ and $H'$ in $\Ss$ such that $\wt{\xi}_t\le H_t-H'_t\le \wt{\zeta}_t$ for all $t \in [0,T]$.}
\ehyp

It was shown in Proposition 3.10 in Dumitrescu et al. \cite{DQS2016} that if $g \in \Hs $, then  the Mokobodzki condition holds if and only if $\II^g \in \Ss$ (or, equivalently, $\JJ^g \in \Ss$). Hence if Assumptions \ref{ass2.1}--\ref{ass2.2} are complemented by Assumption \ref{ass4.2}, then the existence of a unique solution to the DRBSDE \eqref{DRBSDE2} with a fixed generator $g \in \Hs $ is ensured by Lemma \ref{lem4.2}. We may thus state without proof the following result.

\begin{proposition} \label{pro4.3}
Let Assumptions \ref{ass2.1}--\ref{ass2.2} and \ref{xass4.1}--\ref{ass4.2} be satisfied. If the generator $g$ does not depend on $(y,z)$, then there exists a unique solution $(\wt{Y}, Z, \Ku, \Kl)\in \SSAA $ to the doubly reflected BSDE \eqref{DRBSDE2}.
\end{proposition}

We are now ready to demonstrate the existence and uniqueness of a solution to the DRBSDE \eqref{DRBSDE1}.
For a fixed $\beta \ge 0$,  we denote by $\SSb$ the Banach space $(\Sb \times \Lb ,\|\cdot\|_\beta)$ with the norm
\begin{align*}
\|(w,v)\|^2_{\SSb}:=\|w\|^2_{\Sb}+\|v\|^2_{\Lb}.
\end{align*}
We define the mapping $\Psi : \SSb \rightarrow \SSAA$ as follows: for an arbitrary  $(w,v)\in\SSb $ we set $\Psi(w,v):=(Y^{w,v},Z^{w,v},\Kl^{w,v},\Ku^{w,v})$ where $(Y^{w,v}, Z^{w,v}, \Kl^{w,v},\Ku^{w,v})$ is the unique solution to the DRBSDE \eqref{DRBSDE2} with generator $g_t :=g(t,w_t,v_t)$. In addition, we define the projection $\pi: \SSAA \rightarrow \Sb $, which maps $(w,v,l,u) \in \SSAA$ to $\pi (w,v,l,u):=(w,v)\in\Sb$. Finally, the mapping $S: \SSb  \rightarrow \SSb $ is given by $S:=\pi \circ \Psi$. In view of Proposition \ref{pro4.3}, the mapping $\Psi $ (hence also $S$) is well defined under Assumptions \ref{ass2.1}--\ref{ass2.2} and \ref{xass4.1}--\ref{ass4.2}.

\begin{lemma} \label{lem4.3}
If Assumptions \ref{ass2.1}--\ref{ass2.2} and \ref{xass4.1}--\ref{ass4.2} are satisfied, then $S:\SSb \rightarrow \SSb $ is a contraction mapping for some $\beta \ge 0$ and thus $S$ has a unique fixed point $(\wh{Y},\wh{Z})\in\SSb .$
\end{lemma}

\proof
We consider two pairs $(w^1,v^1)$ and $(w^2,v^2)$ in $\SSb$ and, for brevity, we denote their respective images through $S$ by $(Y^1,Z^1)$ and $(Y^2,Z^2)$.  Using the {\it a priori} estimate from Proposition \ref{xpro4.2} and the fact that the generator $g$ is uniformly $m$-Lipschitz continuous with a constant $\wh{L}$, we obtain for arbitrary $\beta\ge 0$ and $\gamma>0$
\begin{align*}
\|Y\|^2_{\Sb}+\|Z\|^2_{\Lb}\le 1158\beta^{-1}\wh{L}^2(C_Q+1)\big(\|w^1-w^2\|^2_{\Sb}+\|v^1-v^2\|^2_{\Lb}\big).
\end{align*}
If we set $\beta=(2\times 1158)\wh{L}^2(C_Q+1)>0$, then we obtain
\begin{align*}
\|S(w^1,v^1)-S(w^2,v^2)\|_\beta^2=\|(Y^1,Z^1)-(Y^2,Z^2)\|_\beta^2\le\frac{1}{2}\|(w^1,v^1)-(w^2,v^2)\|_\beta^2 .
\end{align*}
By the Banach theorem, the mapping $S: \SSb\rightarrow \SSb$ has a unique fixed point $(\wh{Y}, \wh{Z})$.
\endproof

We are now ready to prove the main existence and uniqueness result for the DRBSDE \eqref{DRBSDE1}.

\begin{theorem} \label{the4.1}
If Assumptions \ref{ass2.1}--\ref{ass2.2} and \ref{xass4.1}--\ref{ass4.2} are satisfied, then the doubly reflected BSDE \eqref{DRBSDE1} with
data $(g,\eta,\UU,\xi,\zeta)$ has a unique solution $(Y,Z,\Kl,\Ku)$.
\end{theorem}

\proof
Let the mapping $T:\SSAA \rightarrow \SSAA$ be defined by $T=\Psi \circ \pi$, which means that $T(w,v,l,u)=(Y^{w,v},Z^{w,v},\Kl^{w,v},\Ku^{w,v})$
for all $(w,v,l,u) \in \SSAA$. From Lemma \ref{lem4.3}, we know that the mapping $S$ has a unique fixed point $(\wh{Y},\wh{Z}) \in \Hs \times \Ls$.
From the definitions of $S$ and $T$, it is clear that, for every $(w,v,u,l)\in \SSAA$,
\begin{equation} \label{ttsx}
T(w,v,l,u)=(S(w,v),\Kl^{w,v},\Ku^{w,v}).
\end{equation}
We claim that $(\wh{Y},\wh{Z},\wh{\Kl},\wh{\Ku})$ where $(\wh{\Kl },\wh{\Ku}) :=(\Kl^{\wh{Y},\wh{Z}},\Ku^{\wh{Y},\wh{Z}})$ is a unique fixed point of $T$.
To this end, we observe that the equality $T(\wh{Y},\wh{Z},\wh{\Kl}, \wh{\Ku})=(\wh{Y},\wh{Z},\wh{\Kl},\wh{\Ku} )$ follows from the equality $S(\wh{Y}, \wh{Z})=(\wh{Y}, \wh{Z})$ and \eqref{ttsx} with $(l,u)=(\wh{\Kl},\wh{\Ku})$. To check that a fixed point of $T$ is unique, it suffices to recall that the fixed point of $S$ is unique and use again \eqref{ttsx}. We have thus shown that the DRBSDE \eqref{DRBSDE2} has a unique solution $(\wh{Y},\wh{Z},\wh{\Kl},\wh{\Ku} )$ where $(\wh{\Kl},\wh{\Ku}) :=(\Kl^{\wh{Y},\wh{Z}},\Ku^{\wh{Y},\wh{Z}})$.
\endproof

\subsection{Picard's Iterations for Solutions to Doubly Reflected BSDEs} \label{sec4.3}

To construct the Picard's sequence converging to the solution of a $(Y,Z,\Kl,\Ku)$, we consider the sequence $\{(Y^k, Z^k)\}_{k=1}^\infty$ where $(Y^k, Z^k)=S^k(w,v)$ for some $(w,v) \in \HHb $. The associated sequence $\{(\Kl^k, \Ku^k)\}_{k=1}^\infty$ is obtained by choosing for every $k \in \bnn$ the last two components of the unique solution $(Y^k, Z^k, \Kl^k, \Ku^k)$ to the DRBSDE with generator $g^{k-1}_t :=g(t, Y^{k-1}_t, Z^{k-1}_t)$.  Notice that the existence of a solution is guaranteed under the conditions of Proposition \ref{pro4.3}.  The next result can be demonstrated in a similar way to the proof of Proposition \ref{xpro3.4} and thus we merely sketch the main steps in the proof of Proposition \ref{pro4.4}.

\begin{proposition} \label{pro4.4}
If Assumptions \ref{ass2.1}--\ref{ass2.2} and \ref{xass4.1}--\ref{ass4.2} are valid, then
\begin{equation} \label{converge}
\lim_{k\rightarrow\infty}\|Y^k-Y\|_{\Ss}=0,\quad\lim_{k\rightarrow\infty}\|(\Kl^k-\Kl)-(\Ku^k-\Ku )\|_{\Ss}=0,\quad\lim_{k\rightarrow\infty}\|Z^k-Z\|_{\Ls}=0.
\end{equation}
\end{proposition}

\proof
We first observe that last convergence in \eqref{converge} is clear from Lemma \ref{lem4.3}. If we denote $\wt{Y}^k:=Y^k-Y$, $\wt{Z}^k:=Z^k-Z$,
$\wt{\Kl}^k:=\Kl^k-\Kl$ and $\wt{\Ku}^k:= \Ku^k-\Ku$, then we have
\begin{equation*}
\left\{ \begin{array}[c]{ll}
d\wt{Y}^k_t=-\wt{g}^{k-1}_t\,dQ_t+\wt{Z}^{k,*}_t\,dM_t-d\wt{L}^k_t+d\wt{U}^k_t,\medskip\\\wt{Y}^k_T=0,
\end{array} \right.
\end{equation*}
where $ \wt{g}^{k-1}_t :=g(t,Y^{k-1}_t,Z^{k-1}_t)-g(t,Y_t,Z_t)$. We claim that that there exists a constant $C \ge 0$ such that
\begin{align} \label{showr}
\|\wt{Y}^k\|^2_{\Ss}+\|\wt{\Kl}^k-\wt{\Ku}^k\|^2_{\Ss}\le CJ_k
\end{align}
where
\begin{align*}
J_k:=\|\wt{Y}^{k-1}\|^2_{\Hs}+\|\wt{Y}^k\|^2_{\Hs}+\|\wt{Z}^{k-1}\|^2_{\Ls}+\|\wt{Z}^k\|^2_{\Ls}.
\end{align*}
Assuming \eqref{showr} and using Lemma \ref{lem4.3} and the Banach theorem, we deduce that $\lim_{k\rightarrow\infty}J_k=0$ and thus \eqref{converge} follows. It thus remains to show that \eqref{showr} is valid.

\noindent {\it Step 1.} To show that $\|\wt{Y}^k\|^2_{\Ss}\le CJ_k $ for some constant $C\ge 0$, it suffices to follow Step 1 in Proposition \ref{xpro4.2}.  By applying the It\^o formula to $|\wt{Y}^k_t|^2$ and using the Skorokhod conditions, we get $\int_t^T\wt{Y}^k_{u-}\,d\wt{\Kl}^k_u-\int_t^T\wt{Y}^k_{u-}\,d\wt{\Ku}^k_u\le 0$ for every $t\in [0,T]$, we obtain the inequality
\begin{align*}
|\wt{Y}^k_t|^2&+\int_t^T\|m_u\wt{Z}^k_u\|^2\,dQ_u+\sum_{t<s\le T}|\Delta\wt{Y}^k_u|^2-\int_t^T\|m_u^d(\wt{Z}^k_u)^d\|^2\,dQ_u
\\ &\le 2\int_t^T\wt{Y}^k_u\wt{g}^{k-1}_u\,dQ_u-2 \int_t^T\wt{Y}^k_{u-}\wt{Z}^{k,*}_u\,dM_u.
\end{align*}
The rest follows by repeating the arguments from Step 1 in Proposition \ref{xpro4.2}.

\noindent {\it Step 2.} We claim that $\| \wt{\Kl}^k-\wt{\Ku}^k \|^2_{\Ss}\le C J_k$.
To check this, we observe that
\begin{align*}
\wt{\Kl}^k_t-\wt{\Ku}^k_t=\wt{Y}_0-\wt{Y}_t^k-\int_0^t\wt{g}^{k-1}_u\,dQ_u
+\int_0^t\wt{Z}^{k,\ast}_u\,dM_u.
\end{align*}
and by using the arguments from Step 2 in Proposition \ref{xpro3.4}, we can obtain the desired inequality.
\endproof


\end{document}